\documentclass[12pt]{amsart}

\usepackage{amssymb,amsmath,amsthm}  
\usepackage{graphicx}
\usepackage{bm}
\usepackage{enumerate}
\usepackage[shortlabels]{enumitem}
\usepackage{braket}
\usepackage{amsaddr}
\usepackage{amscd}
\usepackage[dvipsnames]{xcolor}
\usepackage[mathscr]{eucal}
\usepackage{epstopdf}
\usepackage[hidelinks]{hyperref}
\usepackage{subcaption}
\usepackage{cleveref}
\usepackage{relsize}

\hypersetup{
  colorlinks   = true, 
  urlcolor     = PineGreen, 
  linkcolor    = PineGreen, 
  citecolor   = PineGreen 
}

 \allowdisplaybreaks
   
\numberwithin{equation}{section}
\newtheorem{theorem}{Theorem}[section] 
\newtheorem{definition}[theorem]{Definition}
\newtheorem{lemma}[theorem]{Lemma}

\newtheorem{remark}[theorem]{Remark}
\newtheorem{proposition}[theorem]{Proposition}

\newcommand{\tu}{\textup}

\newcommand{\dd}{\displaystyle}
\newcommand{\bfa}[1]{\boldsymbol{#1}} 			%

\newcommand{\e}{\epsilon}
\newcommand{\ddiv}{\textup{div}}     				%
\usepackage[numbers]{natbib}
\usepackage[latin9]{inputenc}
\usepackage[letterpaper]{geometry}
\usepackage{pifont}
\usepackage{float}
\usepackage{color}
\usepackage{adjustbox}
\usepackage{diagbox}
\usepackage{multirow}
\usepackage[makeroom]{cancel}
\usepackage{array}
\usepackage{mathrsfs}
\usepackage{tikz}
\usepackage{comment}
\usepackage{csquotes}

\graphicspath{{./Figures/}}

\usepackage{amsfonts}

\usepackage{color}
\usepackage{adjustbox}
\usepackage{diagbox}
\definecolor{black}{rgb}{0,0,0}

\definecolor{red}{rgb}{1,0,0}

\definecolor{blue}{rgb}{0,0,1}

\numberwithin{equation}{section}

\renewcommand{\div}{\mathop{\rm div}\nolimits}

\newcommand{\dsm}{ \mathrm{d}\sigma}
\newcommand{\dv}{ \mathrm{d}v}
\newcommand{\dx}{ \mathrm{d}x}

\newcommand{\dt}{ \mathrm{d}t}
\newcommand{\ds}{ \mathrm{d}s}
\newcommand{\dmu}{ \mathrm{d}{\mu}}
\newcommand{\dnu}{ \mathrm{d}{\nu}}
\newcommand{\beq}{\begin{equation}}
\newcommand{\eeq}{\end{equation}}
\newcommand{\beqq}{\begin{equation*}}
\newcommand{\eeqq}{\end{equation*}}
\newcommand{\beqas}{\begin{eqnarray*}}
\newcommand{\eeqas}{\end{eqnarray*}}
\newcommand{\bsp}{\begin{split}}
\newcommand{\eesp}{\end{split}}

\makeatletter
\newsavebox{\@brx}
\newcommand{\llangle}[1][]{\savebox{\@brx}{\(\m@th{#1\langle}\)}%
  \mathopen{\copy\@brx\kern-0.5\wd\@brx\usebox{\@brx}}}
\newcommand{\rrangle}[1][]{\savebox{\@brx}{\(\m@th{#1\rangle}\)}%
  \mathclose{\copy\@brx\kern-0.5\wd\@brx\usebox{\@brx}}}

\makeatother

\usepackage[english]{babel}

\title[E\MakeLowercase{ntropic convergence and the linearized limit 
for the} B\MakeLowercase{oltzmann equation with external force}]
{\textsf{\LARGE E\MakeLowercase{ntropic convergence and the linearized limit 
for the} B\MakeLowercase{oltzmann equation with external force}}}

\author[T\MakeLowercase{ina} M\MakeLowercase{ai}]{\Large T\MakeLowercase{ina} M\MakeLowercase{ai$^{*,a,b,c,1}$}}

\begin{document}

\begin{abstract}
 
This paper extends the results regarding entropic convergence and the strong linearized 
limit for the Boltzmann equation (without external force) in [C. David Levermore. Entropic convergence and the linearized limit for the Boltzmann equation. \textit{Communications in Partial Differential Equations}, 18(7-8):1231--1248, 1993] to the case of the Boltzmann equation with external force.  Our starting point is the Boltzmann equation with an external force introduced in [Diogo Ars{\'e}nio and Laure Saint-Raymond. \textit{From the Vlasov--Maxwell--Boltzmann System to Incompressible Viscous Electro-magneto-hydrodynamics}, EMS Press, 2019], we then find new conditions on the force and rigorously prove the maintaining result by Levermore.  More specifically, any sequence of DiPerna-Lions renormalized solutions of the Boltzmann equation with external force are shown to have fluctuations
(about the global Maxwellian equilibrium $M$) that converge entropically (and hence strongly in $L^1$) to the solution of the
linearized Boltzmann equation for any positive time, given that its initial fluctuations
about $M$ converge entropically to the provided $L^2$ initial data of the
linearized equation, where the force can be physically significant.
\end{abstract}
\date{\today}
\maketitle

\noindent \textbf{Keywords.} Boltzmann equation $\cdot$ DiPerna--Lions renormalized solutions $\cdot$ Entropic convergence $\cdot$ Linearized limit $\cdot$ Dissipation equality $\cdot$ External force

\vskip10pt

\noindent \textbf{Mathematics Subject Classification.} 35Q20, 35Q99, 82C22, 82C40, 82D10 

\vfill

\noindent 
$^{*}$Corresponding author: 
\textit{Tina Mai};  
$^a$Institute of Research and Development, Duy Tan University, Da Nang, 550000, Vietnam; $^b$Faculty of Natural Sciences, Duy Tan University, Da Nang, 550000, Vietnam; $^c$ Department of Mathematics and Institute for Scientific
Computation, Texas A\&M University, College Station, TX 77843, USA;\\
\texttt{maitina@duytan.edu.vn}; \texttt{tinamai@tamu.edu}

\vspace{10pt}

\noindent $^1$ This paper builds upon the version previously posted on arXiv:\\ \url{https://arxiv.org/abs/1612.05096}

\newpage


\section{Introduction}\label{sec_intro}

Renormalized solutions (or temporally global weak solutions) of the classical Boltzmann equation are known to exist since the late 1980s, thanks to DiPerna and Lions \cite{LionsCauchy} (at least for the case without the force term).  Thereafter, the convergence of these solutions in the incompressible Navier--Stokes limit was investigated by Bardos, Golse, and Levermore \cite{BGL0, BGL1, BGL}, yielding some degree of success.  In \cite{BGL}, the concept of \textit{entropic convergence} was initially presented and was utilized as a reasonable tool to achieve strong convergence results for fluctuations about 
an absolute equilibrium.  In \cite{Levermore}, making use again the notion of entropic convergence as a major technique, Levermore demonstrated that these fluctuations of DiPerna--Lions renormalized 
solutions converge to a solution of the linearized Boltzmann equation in an $L^2$ space, for the case without 
external force.  It is then natural to think about extending these results in \cite{Levermore} to kinetic equations with external forces 
(thanks to Caflisch) 
and with convex entropy so that the notion of entropic convergence can be applied (thanks to Levermore).   

In this spirit, we start with the Boltzmann equation having an external regular force $F(t,x,v) \in L^1_{\tu{loc}}(\dt \, \dx \, \dv)$ such that \cite[p.~125]{DL16}
\begin{equation}\label{Fconds}
\ddiv_v F = 0\,, \quad F \cdot v = 0\,, \quad \textnormal{and} \quad 
F \in L^1_{\textnormal{loc}}\left(\dt; W^{1,1}_{\textnormal{loc}}( \dx \, \dv)\right)\,,
\end{equation}
whose existence of a (DiPerna--Lions) renormalized 
solution has been justified by Ars\'{e}nio and Saint-Raymond in \cite[Theorem~4.1 on p.~125]{DL16}. 

Our contributions are considering an external force $F(t,x,v) \in L^1_{\tu{loc}}(\dt \, \dx \, \dv)$ such that
\begin{equation}\label{FcondsL2}
\ddiv_v (F) = 0\,, \ F \cdot v = 0\,, \ \tu{and } 
  F \in L^1_{\tu{loc}}\left(\dt \,; W^{1,1}_{\tu{loc}}(\dx \, \dv)\right) \cap L^1_{\tu{loc}}(\dt \,; L^2(M \dx\, \dv)) 
\end{equation}
with physical significance and carefully deriving related analysis results as follows.  Here, the condition
$F \cdot v = 0$ as assumed in \cite[p.~122]{DL16} guarantees that the global Maxwellian $M(v)$ is an equilibrium state of the Boltzmann equation (see \eqref{quadM}).  It is interesting to notice that this assumption $F \cdot v = 0$ also helps simplify the linearized Boltzmann equation to an expected form (Lemma \ref{linearizee}) and ease its analysis, making our paper both original and hold for certain external forces $F(t,x,v)$ that satisfy at least $F \cdot v = 0\,.$  Moreover, the new restriction $F \in L^1_{\tu{loc}}(\dt \,; L^2(M \dx\, \dv))$ in \eqref{FcondsL2} is for showing the relative compactness in Proposition \ref{lbet}.  In this new setting for the Boltzmann equation with such an external force $F$ satisfying \eqref{FcondsL2}, through rigorously proving Proposition \ref{lbet}, we fully demonstrate Proposition \ref{sllt} as our main theorem, which preserves the result 
obtained by Levermore \cite{Levermore} (without external force) regarding the linearization approximation: given any $L^2$ initial data for the linearized (about an equilibrium $M$) Boltzmann equation with external force, any sequence of 
DiPerna--Lions 
renormalized solutions of the Boltzmann equation are shown 
to have fluctuations (about $M$) that converge entropically 
(and therefore strongly in $L^1$) to the solution of the linearized Boltzmann equation 
for all positive times, provided that its initial fluctuations about $M$ converge entropically to the $L^2$ initial data. 

Regarding special significance for physics, we wish to mention two examples about the external forces $F$ meeting the conditions \eqref{Fconds}.  The first example is \textbf{magnetic force} (see \cite[Eq.~(3)]{Bobylev_mag}, \cite[p.~4]{nota_chiara_mag}, \cite[p.~3]{DL16}, \cite[p.~349]{coriolis_taylor}, \cite{hall_coriolis}, for instance):
\[F(t,x, v) = q(v \times B(t,x)) \,,\]
which is the Lorentz force due to a magnetic field $B(t,x)$ only (without electric field) acting on a moving particle of electric charge $q \in \mathbb{R}\,,$ and velocity $v \in  \mathbb{R}^3\,,$ with the vector cross product $\times\,.$ This magnetic force can appear in a plasma or gas of charged particles (as electrons or ions).

The second example is \textbf{Coriolis force} (see \cite{coriolis_be}, \cite[Eq.~(9.49) on p.~348]{coriolis_taylor},  \cite{hall_coriolis}, \cite{euler_exF}, for instance):
\[F(t,x,v) = -2m \, \Omega(t,x) \times v\,,\]
which is an \textit{inertial} force, also known as a fictitious force or pseudo-force (but treated as an \textit{external} force in kinetic theory \cite[p.~1]{CercignaniApp}), experienced by objects moving within a (non-inertial) rotating frame of reference (e.g., Earth's rotation, rotating machinery) with angular velocity vector $\Omega(t,x)\,,$ and we can assume the unit mass $m = 1$ for simplicity.  In applications, such Coriolis force occurs on the Earth's surface \cite{euler_exF} as a result of Earth's rotation \cite{hall_coriolis}, and this force significantly impacts geophysical flows (e.g., atmospheric dynamics and ocean currents) or in rotating machinery with gas particles \cite{coriolis_be}.  For example, with this Coriolis force, for our considering Boltzmann equation, its change is better perceived when one looks at the molecule's trajectory, which is now a diverging spiral for the rotating system rather than a straight line \cite{coriolis_be}.  A remark is that when applied to a fixed medium, the Coriolis force mimics an analogous magnetic field, according to Larmor's theorem (\cite[p.~305]{hall_coriolis}, \cite[p.~21]{cor_mag_si2}, \cite[p.~348]{coriolis_taylor}).   

Our paper's results still hold with more general forces (either external or internal) that satisfy the conditions \eqref{Fconds}.  An example is the internal self-consistent magnetic force in the Vlasov--Maxwell--Boltzmann (without the electric field), but this is mostly an idealized or limited expectation because the electric field will not remain zero in the majority of realistic plasma dynamics \cite{DL16}.

Commonly, many internal (self-consistent, or self-induced, or mean-field) forces (see \cite[Eqs.~(1)--(2)]{caflisch_e}, \cite[p.~196]{VReview}, \cite[Eq.~(2.35) and Eq.~(1.1)]{DL16}, for instance) do not meet the condition $F \cdot v=0$ in \eqref{Fconds}. This condition $F \cdot v=0$ typically means that the force is orthogonal to the particle velocity, like the magnetic or Coriolis forces above, which are usually
external rather than internal.  Therefore, given the condition $F \cdot v=0$ in \eqref{Fconds}, we focus on external forces in this paper.

The condition $\tu{div}_v F = 0$ in \eqref{FcondsL2} as assumed in \cite[p.~122]{DL16} is for verifying the local conservation of mass (Proposition \ref{conserlaws}) and the Boltzmann's $H$ theorem (Proposition \ref{Hthm}).  Note that in the literature, to handle the linearized Boltzmann equation with external forces, for example (also as a more general external force) $F(x,v) = -\nabla_x \phi(x) + \bfa{b}(x,v)$ in \cite{tabataFxv}, the condition $\tu{div}_v F = 0$ is often assumed as well.  Thus, both conditions $F \cdot v =0$ and $\tu{div}_v F = 0$ that we set in \eqref{FcondsL2} are natural restrictions on the external forces in our context.

Section \ref{pre} contains preliminary material regarding the Boltzmann equation with external force 
and its linearization, and having a statement about the DiPerna--Lions renormalized solutions thanks 
to Ars\'{e}nio and Saint-Raymond \cite{DL16}, in our new settings.  We also make some crucial preliminary for the proof of our main theorem (Proposition \ref{sllt}). Since the notion used here is a subset of the one found in the first section of \cite{BGL}, the second section of \cite{Levermore}, 
and Chapter 1 of Part I as well as Chapter 4 of Part II of \cite{DL16}, the reader who is already familiar can skip it.  Section \ref{linlim} reintroduces the concept of entropic convergence (as utilized in \cite{BGL, Levermore, L, DL16}) and encompasses all our key results, especially Proposition \ref{lbet} with new proof and Proposition \ref{sllt} (main theorem) with full justification.  The other proofs of several essential propositions and theorems have been 
shown in \cite{BGL, Levermore, DL16} and are therefore not demonstrated here, but all of these propositions and theorems are expressed in full for the sake of completeness and for considering the external force.  In Section \ref{discuss}, some open problems about extending these results to other settings and 
other kinetic equations are also discussed.  Section \ref{conclude} is for the conclusions.  In Appendix \ref{appendix}, the notion regarding spaces is included.

\section{Notion and preliminaries}\label{pre}

Being motivated by the work \cite{Levermore} of Levermore, in the present paper, we investigate the linearized limit of the Boltzmann equation with 
external force stated by Ars\'{e}nio and Saint-Raymond in \cite[Section~4.1.1]{DL16}.  This equation is of 
the form
\begin{align}
\partial_t f + v \cdot \nabla_x f + F \cdot \nabla_v f & = B(f,f)\,, \label{vb1} \\
f(0,x,v) & = f^{\tu{in}} (x,v) \geq 0\,, \label{vb0}
\end{align} 
with a given external force field $F(t,x,v)$ (in $\mathbb{R}^3$) satisfying at least \cite[p.~122]{DL16}
\[F\,, \ \tu{div}_v F  \in L^1_{\text{loc}}\left(\dt \, \dx \,; L^1(M \dv)\right)\,.
\]
These force field conditions are the very minimum needed to define renormalized solutions of (\ref{vb1}) (see Definition \ref{renorm} below).  Nevertheless, the range of applicability of force fields will be further restricted as follows:
\begin{itemize}
\item we assume $\tu{div}_v F = 0\,,$ mainly for the purpose of verifying the local conservation of mass, and so that the Boltzmann's $H$ theorem holds as in Proposition \ref{Hthm};

\item we assume $F \cdot v = 0\,,$ in order for the global Maxwellian $M(v)$ to be an equilibrium state of (\ref{vb1}).  In our paper, this assumption is also crucial for simplifying the linearized Boltzmann equation to an expected form (Lemma \ref{linearizee}).
\end{itemize}

Here, $t \in [0, \infty)\,, x \in \Omega = 
\mathbb{T}^3\,, v \in \mathbb{R}^3$ (to be regarded as a tangent space to $\Omega$ \cite[p.~77]{VReview}); and the mass density $f(t,x,v) \geq 0$ represents the distribution of particles which, at time $t$, are at position $x$ and have velocity $v$ \cite{Levermore, 4Diogo}.  Note that to maintain simplicity, we shall assume throughout this work that the spatial domain is $\Omega = \mathbb{T}^3$ (not really a subset of $\mathbb{T}^3$) without boundaries \cite{VHproof, Levermore}, therefore avoiding the intricate issue of boundary conditions \cite{Levermore}.  

The Boltzmann collision operator is of the bilinear form \cite{DL16}
\begin{equation}\label{bco1}
B(f,k) = \int_{\mathbb{R}^3} \int_{\mathbb{S}^2} (f' k'_{*} - f k_{*}) \,  b(v-v_{*}, \sigma) \, 
 \dsm \, \dv_{*}\,.
\end{equation}
Thus,
\begin{equation}\label{bcoff}
B(f,f) = \int_{\mathbb{R}^3} \int_{\mathbb{S}^2} (f' f'_{*} - f f_{*}) \,  b(v-v_{*}, \sigma) \, 
 \dsm \, \dv_{*}\,.
\end{equation}
Here,
\[f = f(v)\,, \quad f' = f(v')\,, \quad k_{*} = k(v_{*})\,, \quad k'_{*} = k(v'_{*})\,.\]

The velocities $(v',v'_{*})$ are given by the $\sigma$-\textit{representation} of the collision operator $B$ (\cite[p.~5]{DL16}, \cite[Eq.~(5)]{vrep}, \cite[Eq.~(5)]{VReview}, \cite{chen_natasa}):
\begin{equation}\label{vrep1}
 v' = \frac{v+v_{*}}{2} + \frac{|v - v_{*}|}{2} \sigma \,, \qquad v'_{*} = \frac{v+v_{*}}{2} - \frac{|v - v_{*}|}{2} \sigma \,.
\end{equation}
The parameter $\sigma$ ($\in \mathbb{S}^2$) varies in the unit 2--sphere $\mathbb{S}^2$ (or 2-dimensional sphere) of unit radius \cite[Eq.~(5) on p.~80]{VReview}.  Another useful representation of $B$ is of the form (\cite[Eq.~(2.3)]{Levermore}, \cite[Eq.~(2)]{vrep}, \cite[p.~99]{VReview})
\begin{equation}\label{vrep2}
 v' = v - \langle v - v_{*}, \sigma \rangle \sigma \,, \qquad v'_{*} = v_{*} + \langle v - v_{*}, \sigma \rangle \sigma \,.
\end{equation}
We refer to \cite{vrep} for the link 
between the two representations.  In this work, the 
$\sigma$-representation is employed.  

It is assumed that the collisions are \textit{elastic}: when particle pairs collide, their momentum and energy are conserved, that is,
\begin{equation}\label{quad}
 v + v_{*} = v' + v'_{*}\,, \qquad |v|^2 + |v_{*}|^2 = |v'|^2 + |v'_{*}|^2 \,.
\end{equation}
Thus, one can readily show that the quadruple $(v, v_{*}, v', v'_{*})$ parametrized by 
$\sigma \in \mathbb{S}^2$ gives the collection of all solutions to the system of these four equations (\ref{quad}), where $v \in \mathbb{R}^3$ \cite[p.~5]{DL16}.

In \eqref{quad}, the unprimed velocities ($v\,, v_*$) and primed velocities ($v' \,, v'_*$) \cite[p.~1233]{Levermore} represent possible velocities for two interacting particles either before and after (pre and post (\cite[p.~5]{DL16}, \cite[p.~209]{vrep}, \cite{ShiJin}, \cite{strain})), or after and before (post and pre (\cite[p.~15]{Dilute}, \cite[p.~15]{L}, \cite[p.~79]{VReview}, \cite[p.~630]{mouhot1}, \cite{natasa_gamba_ricardo}, \cite{yan_guo})) elastic binary collision.

Furthermore, we assume that only \textit{uncorrelated} particles are involved in collisions, which means that particles that have already collided are not anticipated to collide again in the future \cite[p.~16]{L}.

The \textit{Boltzmann collision kernel} $b \equiv b(z,\sigma)$ is a nonnegative measurable function, solely relying on $|z|$ and the scalar product $z \cdot \sigma$ (due to the \textit{microreversibility} assumption, that is, microscopic dynamics are time-reversible, or the probability that $(v',v'_*)$ are changed into $(v, v_*)$ in a collision process is the same as the probability that $(v, v_*)$ are changed into $(v', v'_*)$ (see \cite[p.~80]{VReview}, \cite[p.~16]{L})).  It evaluates in some manner the statistical 
repartition of post-collisional velocities, provided the pre-collisional velocities \cite[p.~17]{L}.  Also, it heavily depends on the type of microscopic interactions and is of the following form (\cite[p.~6]{DL16}, \cite[p.~17]{L}):
\[b(z,\sigma) = b \left(|z|, \frac{z}{|z|} \cdot \sigma \right) \geq 0\,.\]
Equivalently \cite[p.~6]{DL16}, with  $z = v - v_*\,,$
\[b(v-v_*,\sigma) = b \left(|v-v_*|, \tu{cos}\theta\right)\,, \quad  \tu{cos}\theta = \frac{v-v_*}{|v-v_*|} \cdot \sigma  \,.\]
This Boltzmann collision kernel $b$ is represented by the specific differential \textit{cross-section} $\Sigma \geq 0$ in the traditional form (\cite[p.~170]{crossphy}, \cite[p.~81]{VReview}):
\begin{equation}\label{cross}
b(v-v_*,\sigma) = |v-v_*| \, \Sigma(v-v_*,\sigma)\,,
\end{equation} 
determined by physics \cite[p.~170]{crossphy}.  Thus, $b$ is often named the \textit{cross-section}, by abuse of language \cite[p.~81]{VReview}.  It is assumed that the support of $b$ is large enough for the Boltzmann's $H$ theorem to hold \cite[p.~1233]{Levermore}.  Without loss of generality, one may limit the deviation angle $\theta$ to the range $\left [0, \dfrac{\pi}{2} \right ]$ (see for instance \cite[p.~120]{VReview}).  Furthermore, we will assume that $b$ satisfies the bounds \cite[Eq.~(2.6)]{Levermore}
\begin{equation}\label{bbounds}
0 \leq b(v - v_{*}\,, \sigma) \leq C \, (1 + |v - v_{*}|^2)\,,
\end{equation}
for some finite constant $C$ independent of $\sigma$.  This requirement is fulfilled by classical Boltzmann 
collision kernels with a small deflection cut-off or Grad's angular cutoff (see \cite[Chapter II]{CercignaniApp}, \cite[p.~1233]{Levermore}), assuming $b(z, \sigma) \in L^1_{\text{loc}} (\mathbb{R}^3 \times \mathbb{S}^2)$, that is, at least locally integrable \cite{L, DL16, LionsCauchy}.  It is as a short-range assumption from the physical perspective \cite[p.~120]{VReview}.

There are well-known facts (\cite[p.~17]{Dilute}, \cite[p.~101 and p.~126]{VReview}) that the \textit{pre-post collisional changes of variables} \cite[p.~101]{VReview} or 
simply \textit{collisional 
symmetries} \cite[p.~8]{DL16}
\begin{align}
& (v, v_{*}\,, \sigma) \mapsto (v', v'_{*}\,, a) \quad \left( \text{with} \quad
\sigma = \frac{v' - v'_{*}}{|v' -v'_{*}|}\,, \quad a = \frac{v-v_{*}}{|v-v_{*}|} \right) \quad \tu{and} \label{csym1} \\
& (v, v_{*}\,, \sigma) \mapsto (v_{*}\,, v, -\sigma) \label{csym2}
\end{align}
are involutive and thus have unit Jacobian determinants.  Here, we have employed $|v'-v'_*|=|v-v_*|$ to \eqref{vrep1} to obtain $\sigma = (v' - v'_{*})/ |v' -v'_{*}|\,;$ and thus, the variable transformation $(v, v_{*}\,, \sigma) \to (v', v'_{*}\,, a)$ in \eqref{csym1} essentially involves the exchange of $(v,v_*)$ and $(v',v'_*)$ \cite[p.~126]{VReview}.  Moreover, as a result of microreversibility, \eqref{csym1} keeps the collision kernel $b$ invariant \cite[p.~101]{VReview}. 

Due to \eqref{vrep1}, the quantity $f'f'_* - ff_*$ in \eqref{bcoff} is invariant under the change of variables $(v, v_{*}\,, \sigma) \to (v_{*}\,, v, -\sigma)$ in \eqref{csym2}.  Hence, if required, one can replace (from the very
beginning) the collision kernel $b$ with its ``symmetrized'' version 
\[\bar{b}(z,\sigma) = b(z,\sigma) + b(z,-\sigma)\,,\] provided $z \cdot \sigma >0$ \cite[p.~120 and p.~124]{VReview}.

Subsequently, making use of the collisional symmetries \eqref{csym1} and \eqref{csym2}, one 
can obtain the following \textit{symmetrized integral} (see \cite[Eq.~(45)]{VReview}, \cite[Eq.~(1.5)]{DL16}, \cite[Eq.~(2.13)]{L}, and \cite[Eq.~(1.11)]{Dilute}, for instance):
\begin{align}\label{sym}
\begin{split}
 & \int_{\mathbb{R}^3} B(f,f) \,  \varphi(v) \, \dv \\
 & = \frac{1}{4} \int_{\mathbb{R}^3 \times \mathbb{R}^3 \times \mathbb{S}^2}  \,  (f' f'_{*} - f f_{*}) \, b(v - v_{*}, \sigma) \,  (\varphi + \varphi_{*} - \varphi' - \varphi'_{*}) \, \dv \, \dv_{*}  \dsm
\,,
\end{split}
\end{align}
where $f(v)$ and $\varphi(v)$ are regular enough such that the integrals in \eqref{sym} exist.  This formula will be referred to as \textit{Boltzmann's weak formulation} \cite[p.~104]{VReview}.  For every appropriate $f(v)\,,$ it follows that the above integral in \eqref{sym} vanishes, that is,  
\begin{equation}\label{B0phi1}
\int_{\mathbb{R}^3} B(f,f) \, \varphi(v) \, \dv=0
\end{equation}
if and only if $\varphi(v)$ on the right-hand side of \eqref{sym} is a \textit{collision invariant} \cite[p.~8]{DL16}, that is, a solution to the functional equation
\begin{equation}\label{invariant}
\varphi(v) + \varphi(v_{*}) = \varphi(v') + \varphi(v'_{*})\,, \qquad 
\forall (v, v_{*}, \sigma) \in \mathbb{R}^3 \times \mathbb{R}^3 \times \mathbb{S}^2\,.
\end{equation}
Obviously, due to \eqref{quad}, any linear combination of $\{1,v_1\,, v_2\,, v_3\,, |v|^2\}$ is a collision invariant.  Furthermore, under very weak assumptions, it can be shown that these are the only collision invariants that could exist \cite[pp.~36--42]{Dilute}.

Let $n(x)$ denotes the outward unit normal to $\partial \Omega$ at $x\,.$ Throughout this work, we assume that any solution $f(t,x,v)$ of the Boltzmann equation (\ref{vb1}) is locally integrable, rapidly decaying in $v\,,$ and $\tu{log} f$ has at most polynomial growth as $|v| \to \infty$ for each $(t,x)\,.$
 
With those assumptions, taking into account the conditions \eqref{Fconds} on the 
given force field $F$ at the beginning of this section, successively multiplying the Boltzmann equation (\ref{vb1}) by the collision invariants ${1,v_1 \, , v_2\,, v_3\,, |v|^2}$, then integrating in velocity (and using the product rule as well as divergence theorem) leads formally to the \textit{local conservation laws}, as follows (see mainly  \cite[Eq.~(1.6)]{DL16} with the force, and also \cite{L, VReview}).
 
\begin{proposition}\label{conserlaws} 
Let $f\equiv f(t,x,v)$ be a solution to the Boltzmann equation \eqref{vb1} which is locally integrable and such that $f$ is rapidly decaying in $v$ for each $(t,x)\,,$ and assume that the given force field $F \equiv F(t,x,v)$ \eqref{Fconds} satisfies $\tu{div}_v F=0$ and $F \cdot v =0\,.$  Then the following local conservation laws hold:
\begin{align}\label{lcl}
 \begin{split}
  & \partial_t \int_{\mathbb{R}^3} f \, \dv + \tu{div}_x \int_{\mathbb{R}^3} v f \, \dv = 0\,, \\
  & \partial_t \int_{\mathbb{R}^3} v f \, \dv + \tu{div}_x \int_{\mathbb{R}^3} (v \otimes v) f \, \dv =\int_{\mathbb{R}^3} F f \, \dv\,, \\
& \partial_t \int_{\mathbb{R}^3} \frac{1}{2}|v|^2 f \, \dv + \tu{div}_x \int_{\mathbb{R}^3} \frac{1}{2}|v|^2 v f \, \dv = 0\,, 
\end{split}
\end{align}
respectively the local conservation of mass, momentum, and energy.
\end{proposition}

\noindent It is well-known that these local conservation laws connect to a description of the gas at the macroscopic level \cite[pp.~8--9]{DL16}.

The symmetries of the collision operator $B$ \eqref{bco1} also yield an additional crucial characteristic of the 
Boltzmann equation.  Without taking into account integrability matters, we plug $\varphi = \text{log} f$ into the symmetrized integral (\ref{sym}) and utilize the logarithm's properties to derive the \textit{entropy production functional} $D(f)$ \cite[p.~4]{VHproof} (or \textit{dissipation of H functional} \cite[p.~4]{VHproof}, or \textit{dissipation of information} according to the conventions discussed in \cite[p.~74]{VReview}, in addition to the terminology \textit{entropy dissipation functional} \cite[p.~74]{VReview}) as \cite[p.~10]{DL16}:
\begin{align}\label{epf}
\begin{split}
 D(f) & =  - \int_{\mathbb{R}^3} B(f,f) \, \text{log} f \, \dv \\
 & =  \frac{1}{4} \int_{\mathbb{R}^3 \times \mathbb{R}^3 \times \mathbb{S}^2} 
(f' f'_{*} - f f_{*}) \, \text{log} \left( \frac{f' f'_{*}}{f f_{*}}\right)\, 
b(v - v_{*}, \sigma)  \, \dv \, \dv_{*} \,  \dsm  \geq 0 \,.
\end{split}
\end{align}

\noindent The so-defined \textit{entropy production} (\cite[p.~251]{H-dissipation}, \cite[p.~55]{VHproof}, \cite{jeremy_carrillo}) (or \textit{entropy dissipation} \cite[p.~10]{DL16}, or \textit{entropy dissipation rate} \cite[p.~1239]{Levermore}, or \textit{$H$-dissipation} \cite[p.~250]{H-dissipation})
\begin{equation}\label{enprod}
\int_{\mathbb{T}^3} D(f)(t,x) \, \dx
\end{equation}
is non-negative, and the functional 
\[\int_0^t \int_{\mathbb{T}^3} D(f)(s,x) \, \dx \, \ds\] 
is thus non-decreasing on $t>0$ \cite[p.~10]{DL16}.

When $B(f,f)=0\,,$ since necessarily $D(f) = 0$ in this case, it is possible to show (\cite[p.~21]{L}, \cite{perthame}) that the minimizers of 
this entropy production functional $D(f)$ are the so-called \textit{global Maxwellian distributions} defined by \cite{Levermore, VReview, DL16} 
\begin{equation}\label{md}
M_{\rho,u,\theta}(v) = \frac{\rho}{(2 \pi \theta)^{\frac{3}{2}}} e^{- \frac{|v - u|^2}{ 2\theta}}\,,
\end{equation}
where $\rho \in \mathbb{R}_{+}\,, u \in \mathbb{R}^3\,, \theta \in \mathbb{R}_{+}\,,$ are respectively 
the constant macroscopic density, bulk velocity, and temperature (which does not depend on $t$ nor $x$) under a suitable selection of units \cite[pp.~10--11]{DL16}.  
The relation $B(f,f) = 0$ exhibits the fact that collisions are no longer the cause of any variation in the density $f\,,$ indicating that the gas has attained statistical \textit{equilibrium} \cite[p.~16]{DL16}.  Thus, a state as \eqref{md} is referred to as a \textit{global Maxwellian} or \textit{global equilibrium} \cite[p.~110]{VReview}, or (\textit{equilibrium}) \textit{absolute Maxwellian}  \cite[p.~1233]{Levermore}, or \textbf{\textit{global Maxwellian equilibrium}} \cite[p.~16]{DL16}.  Such a state \eqref{md} is solely determined by its total mass,
momentum and energy \cite[p.~110]{VReview}.

Let us now introduce the 
\textbf{\textit{Boltzmann's H functional}} (or just the \textbf{\textit{H functional}}, or \textit{information}, or \textit{quantity of information} \cite[p.~74]{VReview}) as 
\begin{equation}\label{hf}
 H(f) := H(f) (t) = \int_{\mathbb{R}^3 \times \mathbb{T}^3} f(t,x,v) \, \text{log} f(t,x,v) \, \dv \, \dx\,,
\end{equation}
which is the \textit{mathematical entropy} (namely, \textbf{\textit{entropy}} \cite[p.~10]{DL16}), convex, and used in this paper.  The \textit{physical entropy} associated with $f$ is defined as \cite{VHthm}
\begin{equation}\label{ef}
 S(f) = - H(f)\,.
\end{equation}

Then, the non-negativeness of the entropy production functional $D(f)$ from \eqref{epf} results in the famous \textbf{\textit{Boltzmann's H theorem}} (\cite[p.~21]{L}, \cite[p.~10]{DL16}, \cite[p.~104]{VReview}), commonly known as the second principle of thermodynamics, which states that the Boltzmann's $H$ functional \eqref{hf} is non-increasing with time, that is (at least formally) a Lyapunov functional of the Boltzmann equation, as follows. 

\begin{proposition}\label{Hthm}
Let $f\equiv f(t,x,v)$ be a well-behaved (smooth) solution of the Boltzmann equation \eqref{vb1} (specifically with finite entropy) which is locally integrable and such that $f$ is rapidly decaying in $v$ and $\tu{log} f$ has at most polynomial growth as $|v| \to \infty$ for each $(t,x)\,,$ and assume that the given force field $F \equiv F(t,x,v)$ \eqref{Fconds} satisfies $\tu{div}_v F=0$ and $F \cdot v =0\,.$  Then the following holds:
\begin{equation}\label{hlyapunov}
\frac{d}{dt} H(f)(t) \leq 0\,.
\end{equation}
Moreover,
\begin{equation}\label{hd}
\frac{d}{dt} H(f)(t) = - \int_{\mathbb{T}^3} D(f)(t, x) \, \dx\,.
\end{equation}
\end{proposition}

\begin{proof}
Formally multiplying the Boltzmann equation \eqref{vb1} by ($\text{log} f +1 $), where $\text{log} f = \text{log}_{e} f$ and recalling $\ddiv_v F=0$ as well as \eqref{sym}--\eqref{invariant}, then integrating 
with respect to velocity and using \eqref{epf} readily lead to
\begin{equation}\label{lya}
 \frac{d}{dt} \int_{\mathbb{R}^3} f \, \text{log} f \, \dv + \ddiv_x \int_{\mathbb{R}^3} 
 v f \, \text{log}f \,  \dv + \int_{\mathbb{R}^3} \ddiv_v (F f \, \text{log} f) \, \dv = - D(f)(t, x) \leq 0\,.
 \end{equation}
On the left-hand side of \eqref{lya}, the last integral of the sum vanishes due to the divergence theorem, the fact that the surface boundary $S_v$ of velocity space $\mathbb{R}^3$ is at infinity, and the hypotheses that $\ddiv_v F=0$ and $f$ is rapidly decaying in $v$ and $\tu{log} f  \leq f$ when $|v| \to \infty$ for each $(t,x)\,.$  (See also the proof in \cite[p.~10]{DL16} for the Vlasov--Maxwell--Boltzmann system with the Lorentz force $F(t,x,v)=q \, (E(t,x) + v \times B(t,x))$ satisfying $\div_v F =0\,.$)    

Next, integrating \eqref{lya} with respect to space, then using the divergence theorem, we get 
\begin{align}\label{to}
\begin{split}
 \frac{d}{dt} H(f)(t) &= \int_{\mathbb{T}^3}\left(\frac{d}{dt} 
 \int_{\mathbb{R}^3} f \,  \text{log} f(t,x,v)\, \dv\right) \, \dx\\
 &= - \int_{\mathbb{T}^3} D(f)(t, x)\, \dx - \int_{\mathbb{T}^3}  \int_{\mathbb{R}^3} 
 \ddiv_x (v f \, \text{log} f(t,x,v)) \, \dv \, \dx\\
 &=  - \int_{\mathbb{T}^3} D(f)(t, x) \, \dx - 
 \int_{\mathbb{R}^3}
 \int_{\mathbb{T}^3}  
 \ddiv_x (v f \, \text{log} f(t,x,v))  \, \dx \, \dv\\
 & = - \int_{\mathbb{T}^3} D(f)(t, x) \, \dx \,. 
 \end{split}
\end{align}
In the third line of \eqref{to}, the second integral on the right-hand side is $0$ due to the periodic boundary condition of the spatial space $\Omega = \mathbb{T}^3$ (no boundaries, $\partial \mathbb{T}^3 = \emptyset$) (see \cite[p.~5, p.~3, and p.~8]{VHproof} where that vanishing integral also happens with smooth bounded $\Omega$ satisfying bounce-back or specular reflection boundary conditions, see \cite{mika_specular} for an example).  Hence, we obtain \eqref{hd}. Finally, the Boltzmann's $H$ theorem \eqref{hlyapunov} holds.
\end{proof}
 
For any particle number density $f \geq 0$ and any global Maxwellian distribution 
$M_{\rho,u,\theta}$ \eqref{md}, 
we define the (global) \textit{relative entropy} of $f$ with respect to
$M_{\rho,u,\theta}$ by \cite[Eq.~(1.17)]{DL16}
\begin{equation}\label{re}
H(f|M_{\rho,u,\theta})(t) = \int_{\mathbb{R}^3 \times \mathbb{T}^3} 
 \left( f \,  \text{log} \left(\frac{f}{M_{\rho,u,\theta}}\right) - f + M_{\rho,u,\theta} \right) (t) \,  \dv \, \dx \geq 0\,,
\end{equation}
which gives a measure of how close $f$ is to such global equilibrium $M_{\rho,u,\theta}$ \cite[p.~1235]{Levermore}. As in \cite[p.~126]{L}, this relative entropy \eqref{re} is a non-negative Lyapunov functional of the Boltzmann equation \eqref{vb1} (see Boltzmann's $H$ theorem in Proposition \ref{Hthm}).  The reason for choosing $H$ as the entropy functional (\ref{re}) is that its integrand is a non-negative strictly convex function of $f$ with a minimum value of $0$ at $f = M_{\rho,u,\theta}$ \cite[p.~1235]{Levermore}.  For simplicity, we denote the relative entropy by $H(f)$, whenever the implied relative Maxwellian distribution $M_{\rho,u,\theta}$ is apparent \cite[p.~11]{DL16}.  

Throughout this work, we will restrict our 
attention to the case of a spatial 
domain without boundary, in particular, a periodic box $\Omega = \mathbb{T}^3 = S^1 \times S^1 \times S^1$ (a torus, which is a compact three-dimensional torus or 3-torus, as the Cartesian product of three circles), thus eliminating all technical obstacles due to boundaries (\cite[p.~1234]{Levermore}, \cite[p.~369]{4Diogo}).  Without loss of generality, we consider linearizations about the dimensionless global normalized Maxwellian equilibrium $M_{1,0,1}(v)$ with $\rho =1\,, u =0\,, \theta =1$  (\cite[p.~16 and p.~122]{DL16}, \cite[p.~1234]{Levermore}) from \eqref{md}, that is,
\begin{equation}\label{m}
 M:=M(v) = \frac{1}{(2 \pi)^{\frac{3}{2}}} e^{-\frac{|v|^2}{2}}\,.
\end{equation}
Note that for all $v \in \mathbb{R}^3\,,$ one has 
\begin{equation}\label{Mbounds}
0 < M(v) \leq \frac{1}{(2 \pi)^{\frac{3}{2}}} < 1\,.
\end{equation}
 
\noindent The condition 
$F \cdot v = 0$ then ensures that this global Maxwellian $M(v)$ \eqref{m} is an equilibrium state of (\ref{vb1}).  
Indeed, plugging $f=M(v)$ into (\ref{vb1}), we obtain
\begin{equation}\label{quadM}
B(M,M) = \partial_t M + v \cdot \nabla_x M + F \cdot \nabla_v M = F \cdot (-v M) 
= (-M)(F \cdot v) = 0\,.
\end{equation}
Finally, the units of length are chosen so that $\mathbb{T}^3$ has a unit volume, leading to the 
normalizations \cite[Eq.~(2.9)]{Levermore}
\begin{equation}\label{norms}
 \int_{\mathbb{S}^2} \, \dsm = 1\,, \qquad \int_{\mathbb{R}^3} M \, \dv = 1\,, \qquad 
 \int_{\mathbb{T}^3} \dx = 1\,.
\end{equation}
Because $M \, \dv$ is a positive unit measure on $\mathbb{R}^3$, we denote by 
$\langle \xi \rangle$ the average over this measure of any integrable function 
$\xi = \xi (v)$ \cite[Eq.~(2.10)]{Levermore}:
\begin{equation}\label{single}
 \langle \xi \rangle = \int_{\mathbb{R}^3} \xi M \, \dv\,.
\end{equation}

From now on, we are interested in the fluctuations of a density $f(t,x,v)$ around $M(v)$ 
defined in (\ref{m}), as a global normalized Maxwellian equilibrium.  It is thus natural to make use of the \textit{relative density} 
$G \equiv G(t,x,v)$ defined by $f=MG$ \cite[p.~1234]{Levermore}.  The initial-value problem \eqref{vb1}--\eqref{vb0} can be recast for $G$ to produce 
\begin{align}
\partial_t G + v \cdot \nabla_x G + F \cdot \nabla_v G & = Q(G,G)\,, \label{vb2} \\
G(0,x,v) & = G^{\text{in}}(x,v) \geq 0\,, \label{vb20}
\end{align}
where we denote the collision operator \cite[p.~123]{DL16} as
\begin{equation}\label{coqb}
Q(G,K) = \frac{1}{M} B(MG,MK)\,.
\end{equation}
Hence, using \eqref{m} and \eqref{quad}, we get $M'M'_* = MM_*\,,$ and thus \eqref{bcoff} becomes 
\begin{equation}\label{bco2}
 Q(G,G) = \int_{\mathbb{R}^3 \times \mathbb{S}^2} (G' G'_{*} - G G_{*}) \, b(v-v_{*}, \sigma) \, 
  M_{*} \, \dv_{*} \, \dsm \,.
\end{equation}
The bound (\ref{bbounds}) on the collision kernel $b$ and the bounds \eqref{Mbounds} ensure that 
\begin{equation}\label{trib}
 \int_{\mathbb{R}^3 \times \mathbb{R}^3 \times \mathbb{S}^2}  b(v - v_{*}, \sigma) \,  
 M_{*} \, \dv_{*} \, M \, \dv \, \dsm < \infty\,.
\end{equation}
Since $\dmu \equiv b(v-v_{*}, \sigma)  M_{*} \, \dv_{*} \, M \, \dv \, \dsm$ is a non-negative finite measure on 
$\mathbb{R}^3 \times \mathbb{R}^3 \times \mathbb{S}^2$ of $(v,v_{*},\sigma)\,,$ we denote by 
$\llangle \Phi \rrangle$ 
the integral over this measure of any integrable function $\Phi = \Phi(v,v_{*},\sigma)$,
\begin{equation}\label{dou}
\llangle \Phi 
 \rrangle = \int_{\mathbb{R}^3 \times \mathbb{R}^3 \times \mathbb{S}^2} \Phi \, \dmu \,.
\end{equation}

If $G$ solves the Boltzmann equation \eqref{vb2}--\eqref{vb20}, then it satisfies the local entropy dissipation 
law (see the derivation in \cite[pp.~671--674]{BGL}):
\begin{align}\label{ledl}
\begin{split}
\partial_t \langle G \, \text{log} G - G  + 1 \rangle + & \ddiv_x \langle 
v \, (G \, \text{log}G - G + 1) \rangle \\
& = -\frac{1}{4} \llangle[\bigg] \text{log} \left(\frac{G' G'_{*}}{G G_{*}} \right) 
(G' G'_{*} - G G_{*}) \rrangle[\bigg] \leq 0 \,.
\end{split}
\end{align}
Indeed, as in the proof of Proposition \ref{Hthm}, mutiplying the Boltzmann equation \eqref{vb2} by $\text{log}G$ and recalling $\ddiv_v F=0$ and $F\cdot v =0$ as well as \eqref{sym}--\eqref{invariant}, then integrating with regard to velocity yields the desired \eqref{ledl}, whose third term on the left-hand side (involving the force $F$) vanishes in a manner similarly to the proof of \eqref{lya}.   Also, note that the local conservation laws in Proposition \ref{conserlaws} holds with $f=MG\,.$ 

Furthermore, integrating \eqref{ledl} over space $\mathbb{T}^3$ and time $[0,t]$ gives the \textit{global entropy equality}
\begin{equation}\label{gein}
 H(G(t)) + \int_0^t R(G(s)) \, \ds = H(G^{\text{in}})\,,
\end{equation}
where $H(G)$ is the relative entropy \cite[Eq.~(2.17)]{Levermore}
\begin{equation}\label{enfunc}
 H(G) := H(G(t))=\int_{\mathbb{T}^3} \langle G \text{log} G - G + 1 \rangle \, \dx \,,
\end{equation}
and $R(G)$ is the entropy production \eqref{enprod} (employing $D(f)$ defined in \eqref{epf}):
\begin{align}\label{endr}
\begin{split}
 R(G) := R(G(s)) &= \int_{\mathbb{T}^3} D(f) (s,x) \, \dx \\
 &= 
\int_{\mathbb{T}^3} \frac{1}{4} \llangle[\bigg] (G'G'_{*} - GG_{*}) \,
 \text{log} \left(\frac{G'G'_{*}}{GG_{*}}\right) \rrangle[\bigg]   \, \dx \geq 0\,.
\end{split}
\end{align}

As mentioned in the definition of the relative entropy \eqref{re}, the reason for choosing $H(G)$ as the functional \eqref{enfunc} is that its integrand is a non-negative strictly convex function of $G$ with a minimum value of $0$ at $G = 1\,.$  
Hence, for any $G\,,$ 
\begin{equation}\label{hg}
H(G) = H(G|1) \geq 0\,, \quad \text{and} \quad H(G) = 0 
\quad \text{iff} \quad G = 1\,.
\end{equation}
This $H(G)$ is known as the relative entropy of $G$ with regard to the absolute equilibrium $G=1$ (that is, in the relative entropy definition \eqref{re}, letting $f=G$ and $M_{\rho,u,\theta} = 1$ will lead to $H(G)$ of the form \eqref{enfunc}); which gives a natural indicator of how close $G$ is to that equilibrium $G = 1$ \cite[p.~1235]{Levermore}.

\subsection{Linearization of the Boltzmann equation with external force}\label{linearizeB}

Linearizing the Boltzmann equation \eqref{vb2}--\eqref{vb20} about $1$, that is, $G = 1 + g$, leads to
\begin{align}
\partial_t g + v \cdot \nabla_x g + F \cdot \nabla_v g + Lg & = 0\,, \label{lbe}\\
g(0,x,v) & = g^{\text{in}}(x,v)\,. \label{lbe0}
\end{align}
Here, the linearized collision operator $L$ (by expanding the operator $Q(G,G)=Q(1+g,1+g)$ 
in (\ref{coqb}), utilizing its bilinearity, and the identity $Q(1,1) =0$) is given by \cite[Eq.~(2.21)]{Levermore} (see also \cite[p.~134]{VReview}, \cite[Eq.~(2.10)]{DL16}, \cite[p.~633]{mouhot1}):
\begin{align}\label{lg}
\begin{split}
 Lg & \equiv   - \frac{1}{M} (B(M,Mg) + B(Mg,M)) = -(Q(1, g) + Q(g,1))\\
 &= \int_{\mathbb{R}^3} \int_{\mathbb{S}^2} (g + g_{*} - g' - g'_{*}) \, b(v - v_{*}, \sigma) \, \dsm \, M_{*} \, \dv_{*} \,
 \end{split}
\end{align}
where
\begin{equation}\label{Qgg}
Q(g,g) = \frac{1}{M}B(Mg,Mg)\,.
\end{equation}

For the Maxwellian equilibrium $M$ defined in (\ref{m}), we have the following result.
\begin{lemma}\label{linearizee}
With the assumption that the given force field $F \equiv F(t,x,v)$ \eqref{Fconds} satisfies $F \cdot v = 0\,,$ we obtain that the linearization of the Boltzmann equation 
 (\ref{vb1}) about $M$ is equivalent to the linearization of 
 the Boltzmann equation (\ref{vb2}) about $1\,,$ and they have the same form \eqref{lbe}.
\end{lemma}
\begin{proof}
Since $f = MG = M(1+g) = M + Mg$, the Boltzmann equation (\ref{vb1}) can be written as
\[\partial_t(M+Mg) + v \cdot \nabla_x(M+Mg) + F \cdot \nabla_v(M+Mg) = B(M+Mg, M+Mg)\,.\]
Equivalently (using the bilinearity of $B$ and the identity $B(M,M)=0$ \eqref{quadM}), we get
\[ M \partial_t g + M v \cdot \nabla_x g + F \cdot (-Mv - Mvg + M \nabla_v g) 
 =  B(M,Mg) + B(Mg,M) + B(Mg,Mg)\nonumber\,.\]
Taking into account the condition $F \cdot v = 0\,,$ the defined collision operator $Q$ (\ref{coqb}), and $Lg$ \eqref{lg}, we reach
\begin{align*}
\partial_t g + v \cdot \nabla_x g + F \cdot \nabla_v g &= \frac{1}{M} (B(M,Mg) + B(Mg,M)) + \frac{1}{M} B(Mg,Mg) \\
&= - Lg + \frac{1}{M} B(Mg,Mg)\,,
\end{align*}
where $B(Mg,Mg)$ can be regarded as very small if $g$ is reasonably near to $0$ \cite[p.~134]{VReview}.  Therefore, linearizing the Boltzmann equation \eqref{vb1} about $M$ also yields the 
linearized equation of the form \eqref{lbe}, and Lemma \ref{linearizee} holds.
\end{proof}

Now, the existence and uniqueness of solution to the linearized Boltzmann equation \eqref{lbe}--\eqref{lbe0} 
are asserted, thanks to the semigroup method discussed in \cite{ScharfSemi} by Scharf, also with external forces $F$ (see a summary in \cite[p.~577]{drange_scharf}).  Indeed, by the arguments in \cite[p.~675]{BGL} and \cite[p.~1236]{Levermore}, the 
operator $L$ has a non-negative definite self-adjoint Friedrichs extension to the Hilbert space 
$L^2(M \dv)\,,$ with the inner product \[\langle g,k \rangle_{L^2(M \dv)} = \langle g k \rangle =  \int_{\mathbb{R}^3} g k M \, \dv \,,\] 
using the notation $\langle \xi \rangle$ in \eqref{single}. Utilizing integration by parts, we get 
\begin{align}\label{skewv}
\begin{split}
\langle (v \cdot \nabla_x)g, k \rangle_{L^2(M \dx \, \dv)} &= \int_{\mathbb{R}^3} \left (\int_{\mathbb{T}^3} (kv) \cdot (\nabla_x g) \, \dx\right) M \, \dv \\
&= 
- \int_{\mathbb{R}^3} \left (\int_{\mathbb{T}^3} (gv) \cdot (\nabla_x k) \, \dx \right ) M \, \dv\\
& = - \langle g, (v \cdot \nabla_x)k \rangle_{L^2(M \dx \, \dv)}\,,
\end{split}
\end{align}
and
\begin{align}\label{skewF}
\begin{split}
\langle (F \cdot \nabla_v)g, k \rangle_{L^2(M \dx \, \dv)} &= \int_{\mathbb{T}^3} \left (\int_{\mathbb{R}^3} (kF) \cdot (\nabla_v g) M \, \dv \right)  \dx \\
&= 
- \int_{\mathbb{T}^3} \left (\int_{\mathbb{R}^3} (gF) \cdot (\nabla_v k) M \, \dv \right )  \dx \\
& = - \langle g, (F \cdot \nabla_v)k \rangle_{L^2(M \dx \, \dv)}\,.
\end{split}
\end{align}
Thus, the operators $v \cdot \nabla_x$ and 
$F \cdot \nabla_v$ are skew-adjoint on $L^2(M \dx \, \dv)$, and 
\[-v \cdot \nabla_x - F \cdot \nabla_v - L\]
is a closable dissipative operator on $L^2(M \dx \, \dv)$ that generates a strongly continuous contraction 
semigroup.  Hence, for every $g^{\text{in}} \in L^2(M \dx \, \dv)$, there exists a unique 
\[g \in C\left([0, \infty)\,; L^2(M \dx \, \dv)\right)\]
that solves the linearized Boltzmann equation \eqref{lbe}--\eqref{lbe0}.  For each such solution $g(s,x,v)\,,$ the weak formulation of \eqref{lbe}--\eqref{lbe0} holds, as follows:
\begin{align}\label{wle}
\begin{split}
 \int_{\mathbb{T}^3} \langle g(t_2) \varphi \rangle \, \dx - 
 \int_{\mathbb{T}^3} \langle g(t_1) \varphi \rangle \, \dx & -
 \int_{t_1}^{t_2} \int_{\mathbb{T}^3} 
 \langle g(v \cdot \nabla_x \varphi + F \cdot \nabla_v \varphi)  \rangle \, \dx \, \ds \\
 & + \int_{t_1}^{t_2} \int_{\mathbb{T}^3} \langle (Lg) \varphi \rangle \, \dx \, \ds = 0\,, 
\end{split}
\end{align}
for every test function $\varphi(x,v) \in W^{1,2}\left(M \dv\,; C^1(\mathbb{T}^3)\right)$ (see \cite[Eq.~(2.22)]{Levermore}), and for every $0 \leq t_1 < t_2 < \infty\,.$  The solution $g$ also satisfies the 
dissipation equality \cite[Eq.~(2.23)]{Levermore}: 
\begin{equation}\label{disse}
\int_{\mathbb{T}^3} \frac{1}{2} \langle g^2(t) \rangle \, \dx +
\int_0^t \int_{\mathbb{T}^3} \frac{1}{4} \llangle  q^2 \rrangle \, \dx \, \ds 
= \int_{\mathbb{T}^3} \frac{1}{2} \langle (g^{\text{in}})^2\rangle \, \dx\,,
\end{equation}
where $q = (g' + g'_{*} - g - g_{*})$ satisfies the classical identity 
$\langle g (Lg) \rangle = \frac{1}{4} \llangle q^2 \rrangle$ (see \cite[Eq.~(1.44) on p.~675]{BGL}), with the notation $\llangle \cdot 
 \rrangle$ defined in \eqref{dou}.  Indeed, to prove \eqref{disse}, choosing the test function $\varphi = g(s,x,v)\,,$ 
 and multiplying \eqref{lbe} by this test function, where $s \in [0,t]\,,$ we obtain the weak formula of \eqref{lbe}, as follows:
\begin{equation}\label{dissg}
 \int_{0}^{t}  \int_{\mathbb{T}^3} \langle g \partial_t g  \rangle \, \dx  -
 \int_{0}^{t} \int_{\mathbb{T}^3} 
 \langle g(v \cdot \nabla_x g + F \cdot \nabla_v g)  \rangle \, \dx \, \ds + \int_{0}^{t} \int_{\mathbb{T}^3} \langle (Lg) g \rangle \, \dx \, \ds = 0\,.
 \end{equation}
 Recall that the operator $L$ is self-adjoint on $L^2(M \dv)\,,$ and the operators $v \cdot \nabla_x$ and 
$F \cdot \nabla_v$ are skew-adjoint on $L^2(M\dx \, \dv)$ as \eqref{skewv}--\eqref{skewF}, so the second term on the left-hand side of \eqref{dissg} vanishes.  Hence, \eqref{disse} holds as desired. 
 

\bigskip

For later use, recall that the prefixes $w\textnormal{-}$ or $w^{*}\textnormal{-}$ express that a given space is equipped with its weak or 
weak-$*$ topology, respectively.  About the notation regarding spaces, see Appendix \ref{appendix}. 

\subsection{Renormalized solutions of the Boltzmann equation with external force}\label{secRS}

DiPerna and Lions \cite{LionsCauchy} proved the existence of a renormalized solution (also known as temporally global weak solution \cite[p.~1236]{Levermore}) to the Boltzmann equation without external force (that is, the case $F=0$ in \eqref{vb1} and \eqref{vb2}) 
over the spatial domain $\mathbb{R}^3\,,$ given arbitrarily large initial data 
$G(0,x,v) = G^{\text{in}}(x,v) \geq 0$ for \eqref{vb2}--\eqref{vb20}, satisfying natural physical bounds.  With slight modifications, their theory 
can be extended to the periodic box $\mathbb{T}^3\,.$  The existence of DiPerna--Lions renormalized solutions for the Vlasov--Boltzmann systems with self-consistent Vlasov-type force field (derived from a self-induced potential) has been established in \cite{Lions3} and additionally studied then (see \cite{mischler_vlasov,Diogo13Sol}, for instance).  

For the 
Boltzmann equation with external force \eqref{vb2}--\eqref{vb20}, in the recent work \cite{DL16}, 
Ars\'{e}nio and Saint-Raymond have justified the existence of such a renormalized solution, following the construction in \cite{LionsCauchy}.  Toward investigating the linearized limit, we state again these 
results from \cite[Section 4.1.1]{DL16} in our context, with the benefit of the notations and results from \cite{Levermore} by Levermore.  Recall that $F$ is a given force field satisfying \eqref{Fconds}, that is, 
\[F(t,x,v) \in L^1_{\tu{loc}}(\dt \, \dx \, \dv )\]
such that  
\begin{equation*}
\ddiv_v F = 0, \qquad F \cdot v = 0 \qquad \textnormal{and} \qquad 
  F \in L^1_{\textnormal{loc}}\left(\dt\,; W^{1,1}_{\tu{loc}}(\dx \, \dv)\right)\,. 
\end{equation*}
The following definition of renormalized solution is reformulated from \cite{Diogo13Sol, 4Diogo, nader, DL16, LionsCauchy, DL, BGL, Levermore}.  

\begin{definition}\label{renorm} 
Consider the Boltzmann equation \eqref{vb2}--\eqref{vb20} with the external force $F \equiv F(t,x,v)$ satisfying \eqref{Fconds}.  A nonlinearity $\beta(z) \in C^2([0,\infty)\,;\mathbb{R})$ is called an admissible 
\textbf{renormalization} and $N(z) \in C([0,\infty)\,; (0, \infty))$ is named an 
admissible (DiPerna--Lions) \textbf{normalization} if they 
satisfy, for some finite constant $C>0\,,$
\begin{equation}\label{adre}
|\beta'(z)| 
= \frac{1}{|N(z)|} 
\leq 
\frac{C}{1+z} \tu{ and } \beta''(z) \leq 0\,,  
\tu{ for all } z \geq 0 \,,
\end{equation}
where in our paper, we consider the case \cite[p.~1236]{Levermore} 
\[0 < \beta'(z) = \frac{1}{N(z)} \leq 
\frac{C}{1+z} \tu{ for every } z \geq 0 \,.\]

A relative density function $G(t,x,v) \geq 0$, 
where $(t,x,v) \in [0, \infty) \times 
\mathbb{T}^3 \times \mathbb{R}^3\,,$ such that 
\begin{equation}\label{gconditions}
G \in C \left([0,\infty)\,;w\tu{-}L^1(M \dx \, \dv)\right) 
\cap L^{\infty}\left([0,\infty), \dt\,; L^1\left((1+|v|^2) M \dx \, \dv \right)\right)
\end{equation}
is a \textbf{renormalized solution of the Boltzmann equation} \eqref{vb2}--\eqref{vb20} if it 
solves the renormalized Boltzmann equation \cite{DL16,Levermore}
\begin{align}
\partial_t \beta(G) + v \cdot \nabla_x \beta(G) + F \cdot \nabla_v \beta(G) & = 
 \frac{1}{N(G)} Q(G,G)\,, \label{rbe} \\
 G(0,x,v) & = G^{\textnormal{in}}(x,v) \geq 0\,, \label{rbe0}
\end{align}
in the sense of distributions for any admissible renormalization (correspondingly, normalization), 
and satisfies the entropy 
inequality \cite[p.~123]{DL16}, for all $t \geq 0\,,$
\begin{equation}\label{eif}
H(f(t)) + \int_0^t \int_{\mathbb{T}^3} D(f(s)) \, \dx \, \ds \leq H(f^{\textnormal{in}}) < \infty\,,
\end{equation}
where $f^{\textnormal{in}} = MG^{\textnormal{in}}$ is initial value of $f = MG$ and the relative entropy 
$H(f) = H(f|M)$ is defined in (\ref{re}), 
while the entropy production functional (entropy dissipation functional) $D(f)$ is defined in \eqref{epf}.  In this 
definition, equivalently, $f=MG$ is called a renormalized solution \cite{DL16} of the Boltzmann equation \eqref{vb1}--\eqref{vb0}.
\end{definition}

Notice from \cite[p.~124]{DL16}, \cite[Eq.~(2.25)]{Levermore}, and \cite[Eq.~(5.8)]{4Diogo}, that the renormalized collision operator $\beta'(G) Q(G,G)$ is well-defined in 
\[L^1_{\text{loc}}\left(\dt\,; L^1(M \dx \, \dv)\right)\]
for any admissible renormalization (correspondingly, 
normalization), for any function $G$ in (\ref{gconditions}), and any integrable collision kernel 
$b(z, \sigma) \in L^1_{\text{loc}}(\mathbb{R}^3 \times \mathbb{S}^2)$ satisfying the so-called 
DiPerna--Lions assumption \cite[p.~124]{DL16}
\begin{equation}\label{dlassump}
 \lim_{|v| \to \infty} \frac{1}{|v|^2} \int_{K \times \mathbb{S}^2} b(v-v_{*},\sigma) \, \dv_{*} \, \dsm =0\,,
\end{equation}
for any compact set subset $K \subset \mathbb{R}^3\,.$

Indeed, we first consider \textit{non-negative renormalizations} $\beta(z)$ that \cite[p.~1236]{Levermore} 
\[0 < \beta'(z) \leq \frac{C}{1+z}\,.\]
The collision operator $Q(G,G)$ \eqref{bco2} can be decomposed into gain (source) and loss (sink) parts respectively (\cite[p.~81, p.~154]{VReview}, \cite[Eq.~(2.26)]{Levermore}, and \cite[p.~6]{DL16}): 
\begin{equation}\label{deq}
 \begin{split}
  Q^{+}(G,G) = \int_{\mathbb{R}^3 \times \mathbb{S}^2} G' G'_{*} \, b(v-v_{*}, \sigma)  \, M_{*} \, \dv_{*} \, \dsm \,,\\
  Q^{-}(G,G) = \int_{\mathbb{R}^3 \times \mathbb{S}^2} G G_{*} \, b(v-v_{*}, \sigma)  \, M_{*} \, \dv_{*} \, \dsm\,,
\end{split}
\end{equation}
and clearly these gain part $Q^{+}$ and loss part $Q^{-}$ are positive operators \cite{vasseur_mischler}.
From the hypothesis (\ref{dlassump}), it is possible to show directly 
that (as \cite[p.~124]{DL16} and see \cite{4Diogo}, for instance, for more details)
\[\lim_{|v| \to \infty} \frac{1}{|v|^2} 
\int_{\mathbb{R}^3 \times \mathbb{S}^2} b(v-v_{*},\sigma) \, M_{*} \, \dv_{*} \, \dsm =0\,.\]
Thus, the renormalized loss part clearly satisfies \cite{Levermore, 4Diogo}
\begin{align}\label{lossart}
\begin{split}
\beta'(G) \, Q^{-}(G,G) = &\beta'(G) \, G  \int_{\mathbb{R}^3 \times \mathbb{S}^2}  G_{*} \, b(v-v_{*}, \sigma)  \, M_{*} \, \dv_{*} \, \dsm\,.\\
& \in L^{\infty}(\dt \,; L^1(M \dx \, \dv))\,.
\end{split}
\end{align}
Whereas, the renormalized gain part
\begin{align}\label{gainp}
\beta'(G) \, Q^{+}(G,G) = \beta'(G) \int_{\mathbb{R}^3 \times \mathbb{S}^2} G' G'_{*} \, b(v-v_{*}, \sigma)  \, M_{*} \, \dv_{*} \, \dsm
\end{align}
is well-defined in the space $L^1_{\text{loc}}
\left(\dt\,; L^1(M \dx \, \dv)\right)$ \cite{DL16, Levermore, 4Diogo} by estimating a priori \cite[p.~411]{4Diogo} from the renormalized Boltzmann equation (\ref{rbe}) since it is the only remaining unestimated term and it is non-negative \cite{DL16}.  These controls are easily extended to \textit{signed renormalizations} $\beta(z)$ 
satisfying \eqref{adre} with 
\[|\beta'(z)| \leq \frac{C}{1+z}\,,\]
as the Boltzmann equation (\ref{rbe}) is linear with respect to renormalizations so that $\beta'(z)$ can be decomposed into positive and negative parts \cite[p.~124]{DL16}.

In \cite[Eq.~(5.11)]{4Diogo} and \cite[Eq.~(4.8)]{DL16}, an alternative method based on 
entropy dissipation control, which was 
originally performed in \cite{LionsCauchy}, is also discussed to claim that 
the renormalized gain part $\beta'(G) \, Q^{+}(G,G)$ is well-defined in $L^1_{\text{loc}}
\left(\dt\,; L^1(M \dx \, \dv)\right)\,.$ 

Last, note that it is possible to extend the definition of the \textit{renormalized collision operator} 
\begin{equation}\label{renormQ}
\beta'(G) \,Q(G,G)
\end{equation}
to all admissible renormalizations $\beta(z)$ (that is, all $\beta(z)$ satisfying \eqref{adre})  \cite[pp.~124--125]{DL16}.
   
Therefore, saying $G$ is a \textit{weak solution} of the renormalized Boltzmann equation \eqref{rbe} means that it is initially equal to $G^{\text{in}}$ in \eqref{rbe0}, and it should satisfy \eqref{rbe} in the manner that for any non-negative test function
$\varphi (x,v) \in W^{1,\infty}\left(M \dv\,;C^1(\mathbb{T}^3)\right)$ (see \cite[Eq.~(2.27)]{Levermore}, \cite[p.~125]{DL16}), and for every $0 \leq t_1 < t_2 < \infty\,,$ 
\begin{align}\label{wsol}
 \begin{split}
\int_{\mathbb{T}^3} \langle \beta(G(t_2)) \, \varphi \rangle \, \dx 
  -\int_{\mathbb{T}^3} \langle \beta(G(t_1)) \, \varphi \rangle \, \dx
  & - \int_{t_1}^{t_2} \int_{\mathbb{T}^3} \langle \beta(G) \, 
  (v \cdot \nabla_x \varphi + F \cdot \nabla_v \varphi)  \rangle \, \dx \, \dt\\
  & = \int_{t_1}^{t_2} \int_{\mathbb{T}^3} 
  \left \langle \frac{1}{N(G)} Q(G,G) \, \varphi \right \rangle \, \dx \, \dt\,,
 \end{split}
\end{align}
thanks to  the hypothesis $\beta \in C^2([0,\infty)\,;\mathbb{R})$ in the Definition \ref{renorm}.  If integrating over $t\in [0,\infty)\,,$ then one may additionally choose any test functions $\rho(t) \in C_c^{\infty}([0,\infty))$ together with such $\varphi (x,v) \in W^{1,\infty}\left(M \dv\,;C^1(\mathbb{T}^3)\right)$ (see \cite[p.~125]{DL16} and \cite[p.~414]{4Diogo}, also for the corresponding weak formulation). Toward later use, note that
\begin{equation}\label{phisp}
 W^{1,\infty}\left(M \dv\,;C^1(\mathbb{T}^3)\right) \subset W^{1,\infty}(M \dx \, \dv) \subset L^{\infty}(M \dx \, \dv)\,.   
\end{equation}

\bigskip

The following theorem is an extension of \cite[Theorem 4.1 on pp.~125--126]{DL16}, which is a modern formulation of the 
existence result found in \cite{LionsCauchy}.

\begin{theorem}[\textbf{\cite[Theorem 4.1]{DL16}, \cite{LionsCauchy, DL,4Diogo}}]\label{rthm}
Let $b(z, \sigma) \in L^1_{\tu{loc}}(\mathbb{R}^3 \times \mathbb{S}^2)$ be a locally integrable collision kernel satisfying the DiPerna--Lions assumption (\ref{dlassump}), and let $F$ be a given force field satisfying \eqref{Fconds}, that is, 
\[F(t,x,v) \in L^1_{\tu{loc}}(\dt \, \dx \, \dv )\]
such that 
\begin{equation*}
\ddiv_v F = 0, \qquad F \cdot v = 0 \qquad \textnormal{and} \qquad 
  F \in L^1_{\textnormal{loc}}\left(\dt\,; W^{1,1}_{\tu{loc}}(\dx \, \dv)\right)\,.
\end{equation*}
 
Then, for any initial condition $f^{\textnormal{in}} = M G^{\textnormal{in}} \in L^1\left((1+|v|^2) M \dx \, \dv \right)$ such that $f^{\textnormal{in}} = M G^{\tu{in}} \geq 0$ and 
 \[H(f^{\tu{in}}) = H(f^{\tu{in}} | M)= \int_{\mathbb{T}^3 \times \mathbb{R}^3} 
 (G^{\tu{in}} \, \tu{log} G^{\textnormal{in}} - G^{\textnormal{in}} + 1)  M \dx \, \dv
 < \infty\,,\]
 there exists a renormalized solution $f(t,x,v) = M(v) \, G(t,x,v)$ to the Boltzmann equation \eqref{vb1}--\eqref{vb0} (equivalently, a renormalized solution $G(t,x,v)$ to \eqref{vb2}--\eqref{vb20}).  
 Moreover, it satisfies the local conservation law of mass 
\begin{equation}\label{lcm2}
 \partial_t \langle G \rangle + \ddiv_x \langle v G \rangle = 0 \,,
 \end{equation}
 and the global entropy inequality, for any time $t \geq 0\,,$
\begin{equation}\label{gei2}
  H(G(t)) + \int_0^t R(G(s)) \, \ds \leq H(G^{\textnormal{in}})\,.
 \end{equation}
Here, we have 
$H(G) = H(G|1)=H(f) = H(f | M)\,,$ where the relative entropy $H(f):=H(f|M)$ is defined in \eqref{re}, and thus the functional $H(G):=H(G|1)$ has the form (\ref{enfunc}); 
while the entropy production (entropy dissipation rate) $R(G)$ is defined in (\ref{endr}).
\end{theorem}

Some comments on the temporal continuity of the renormalized solutions $G(t)$ (that is, $G\in  \left([0,\infty)\,;w\tu{-}L^1(M \dx \, \dv)\right)$ in \eqref{gconditions}) and the temporal continuity of the renormalizations $\beta(G)$ were discussed in \cite[pp.~415--417]{4Diogo}. 

\begin{remark}
The DiPerna--Lions theory, however, does not assert the local conservation of momentum, the global conservation of energy, nor the global entropy equality (\ref{gein}); neither does it claim that the solution is unique \cite[p.~1238]{Levermore}. In the case of long-range interactions (that is, when the collision kernel is non-integrable), there 
 exist renormalized solutions with a \textit{defect measure} (in the spirit of the construction 
 by Alexandre and Villani \cite{V_long} for the Boltzmann equation), which has been addressed by Ars\'{e}nio and Saint-Raymond 
 \cite{Diogo13Sol} for the Vlasov--Maxwell--Boltzmann system.
 \end{remark}

This Theorem \ref{rthm} can be proved by using the standard procedure in the analysis of weak solutions of 
partial differential equations (see the work of DiPerna and Lions \cite{LionsCauchy, DL} and Lions \cite{Lions1, Lions2, Lions3}); that is to say, solving an approximate truncated equation, determining uniform a priori estimates as well as the weak compactness of the approximate solutions, and ultimately recovering the original equation by passing to the limit (via demonstrating the nonlinear terms' weak stability) \cite[p.~126]{DL16}.  
In many cases, this process reduces to examining the essential weak stability of solutions.  Also, recently, this method of proof in \cite{Lions1, Lions2, Lions3} has been explained and justified in \cite{DL16} by studying the weak stability of renormalized solutions for the Boltzmann 
equation \eqref{vb2}--\eqref{vb20}, equivalently, the weak stability of weak solutions of the renormalized 
equation \eqref{rbe}--\eqref{rbe0}.  We summarize here the crucial steps in the proof of this Theorem 
\ref{rthm} as in \cite[pp.~126--131]{DL16}, following Lions' strategy \cite{Lions1, Lions2, Lions3}, which depends on velocity averaging lemmas, substantial renormalization techniques, and above on, the compactifying (and in certain situations, regularizing) effect of the collision operator's gain term $Q^{+}(G,G)\,.$

\vspace{12pt}
\noindent \textbf{Step 1:}  We consider a 
sequence $\{G_j\}_{j \in \mathbb{N}}$ of actual renormalized solutions to the Boltzmann equation \eqref{vb2}, whose 
initial data $\{G^{\text{in}}_j\}_{j \in \mathbb{N}}$ \eqref{vb20} converges weakly (say, at least 
in $L^1_{\text{loc}}$) when $j \to \infty$ to $G^{\text{in}}$.  Furthermore, we assume 
that the initial data satisfies the following strong entropic convergence
\[\lim_{j \to \infty} H(G_j^{\text{in}}) = H(G^{\text{in}}),\]
so that the entropy inequality is uniformly satisfied
\begin{equation}\label{gei3}
H(G_j(t)) + \int_0^t R(G_j(s)) \, \ds \leq H(G_j^{\text{in}})\,,
\end{equation}
where the relative entropy $H(G_j(t))$ has the form \eqref{enfunc}.

Note from \eqref{gei3} that a uniform bound on the entropies $H(G_j(t))\,,$ through directly using the Young inequality \cite[Eq.~(B.3)]{DL16}, leads to a uniform bound on $G_j(t,x,v)$ in $L_{\tu{loc}}^{\infty}\left(\dt\,; L^1\left((1+|v|^2) M \dx \, \dv \right)\right)\,.$  Also, utilizing the Young inequality \cite[Eq.~(B.3)]{DL16} again together with the Dunford-Pettis compactness criterion (see \cite{Royden} and \cite[Section 5.1]{DL16}), one can show that the $G_j$'s are weakly relatively compact in $L^1_{\text{loc}}\left(\dt \, ; L^1(M \dx\, \dv)\right)\,.$  Thus, extracting a subsequence if necessary, we may assume that the sequence 
$\{G_j\}_{j \in \mathbb{N}}$ converges weakly (as $j \to \infty$) to some $G$ in 
$L^1_{\text{loc}}\left(\dt \,; L^1(M \dx \, \dv)\right)$ (see \cite[Theorem~4.30 on p.~115, and Exercise 3 on p.~468]{brezis}).

One now uses
the idea of 
renormalized solutions to obtain $L^1$ bounds on nonlinear terms 
\[\dd \frac{M \, Q^{\pm}(G_j,G_j)}{1+\delta M G_j *_v \left(\int_{\mathbb{S}^2} b(\cdot , \sigma)  \, \dsm \right)} \qquad \text{and} \qquad 
\beta'(G_j) \, M \,  Q^{\pm}(G_j,G_j)\,, \]
for any $\delta > 0$ and any admissible nonlinearity $\beta(z) \in C^1([0,\infty);\mathbb{R})\,.$  These operators even have 
``weakly compact in $L^1$'' estimates employing the entropy production functional $D$ \eqref{epf}.

\vspace{12pt}
\noindent \textbf{Step 2:} Such estimates are shown to demonstrate that various \textit{velocity averages} 
of $G_j\,,$ $Q^{\pm}(G_j, G_j)$ are relatively compact in $L^1_{\tu{loc}}(\dt \ \dx)$ or almost everywhere, and converge to 
the averages of $G\,,$ $Q^{\pm}(G, G)\,.$  This step relies on a standard use of the velocity averaging 
theory for linear transport equation (as treated in \cite{Dvaverages}, for instance, with $F \cdot \nabla_v \beta(G_j) = 
\ddiv_v(F \,  \beta(G_j))$ as a source term \cite[p.~128]{DL16}).

At this stage, using the convexity methods from \cite{DL}, 
by assuming $\ddiv_v F = 0$ and $F \cdot v=0$, passing to the limit in (\ref{gei3}), one obtains the entropy inequality 
(\ref{gei2}), that is, it follows from formal estimates on the Boltzmann equation (\ref{vb2}), even when 
$F \neq 0\,.$ 

\vspace{12pt}
\noindent \textbf{Step 3:} Employing the information from the earlier steps to reach 
the limit through a somewhat complicated process, which involves interpreting the Boltzmann equation 
(\ref{vb2}) along almost particle paths, and thus recover the original equation \eqref{rbe}. 

In this step, notice that the hypothesis on $F$ that
$F \in L^1_{\textnormal{loc}}\left(\dt; W^{1,1}_{\tu{loc}}(\dx \, \dv)\right)$
is included in the hypotheses of Theorem \ref{rthm} so that the vector field 
$(v, F(t,x,v)) \in \mathbb{R}^6$ satisfies the conditions imposed on 
the transport equation \cite[Eq.~(4.16) on p.~129]{DL16} (see also \cite[Theorem II.1 on p.~516]{DiPerna1989}). \qed

\section{The linearized limit with external force}\label{linlim}
This section of the case with external force will extend the work  of Levermore \cite{Levermore} (without external force).

For small $\epsilon >0\,,$ let $G_{\epsilon} \geq 0$ be a family of DiPerna--Lions renormalized solutions to the scaled Boltzmann 
initial-value problem \eqref{vb2}--\eqref{vb20} so that the initial data $G_{\epsilon}^{\text{in}}$ satisfies 
the \textbf{relative entropy bound} (\cite[p.~158]{DL16}, \cite[Eq.~(3.1)]{Levermore}, \cite[Eq.~(8.1)]{4Diogo}) 
\begin{equation}\label{eb}
H(G_{\epsilon}^{\text{in}}) \leq C^{\text{in}} \epsilon^2\,,
\end{equation}
for some fixed $C^{\text{in}} > 0$ not depending on $\e\,,$ with the relative entropy $H$ specified in \eqref{enfunc}.  Consider the sequence of fluctuations $g_{\epsilon}$ defined by the relation
\begin{equation}\label{ge}
 G_{\epsilon} = 1 + \epsilon g_{\epsilon}\,,
\end{equation}
where $g_{\epsilon} \in L^{\infty} \left( \dt\,; L^1 \left( \left(1+ |v|^2\right) M \dx \, \dv \right) \right)$ \cite[Eq.~(7.1)]{4Diogo}.
The entropy bound \eqref{eb} and the DiPerna--Lions entropy inequality \eqref{gei2} are compatible with this order of fluctuation about the equilibrium $G=1\,.$  In particular, 
in \cite{Levermore}, it is shown that the so-defined sequence $g_\epsilon$ is of order one.

Now, the DiPerna--Lions normalization is chosen to be in the form
\begin{equation}\label{normalization}
 N_{\epsilon} = N(G_{\epsilon}) = \frac{2}{3} + \frac{1}{3} G_{\epsilon} 
 = 1 + \frac{1}{3} \epsilon g_{\epsilon}\,.
\end{equation}
One reason to pick this form is that formally, $N_{\epsilon} \to 1$ when $\epsilon \to 0 $; hence conveniently, the normalizing factor 
will disappear from every algebraic expression taken into consideration in such limit.  Another reason 
is the simplification of the details included in some later estimations.  Our primary findings, of course, are not affected by this specific normalization choice.  Provided this selection \eqref{normalization}, we let 
\begin{equation}\label{gamma}
 \gamma_{\epsilon} = \frac{1}{\epsilon} \beta(G_{\epsilon}) = \frac{3}{\epsilon} \, \tu{log} 
 \left(1+ \frac{1}{3} \epsilon g_{\epsilon} \right).
\end{equation}
The renormalized Boltzmann equation \eqref{rbe} thus becomes 
\begin{equation}\label{gammab}
 \partial_t \gamma_{\epsilon} + (v \cdot \nabla_x + F \cdot \nabla_v) \gamma_{\epsilon} 
 =  \frac{1}{\epsilon} \frac{Q(G_{\epsilon}, G_{\epsilon})}{N_{\epsilon}}\,.
\end{equation}
Note from here that $\gamma_{\epsilon}$ should be viewed as the normalized form of the fluctuations $g_{\epsilon}$  because it formally behaves like $g_{\epsilon}$ for small $\epsilon >0\,.$

The first goal is to describe the characteristics of the limit of the fluctuations $g_{\epsilon}\,.$  For this purpose, one can start to obtain the required a priori estimations by combining the entropy inequality \eqref{gei2} and the entropy bound \eqref{eb} for the initial data, as follows:
\begin{equation}\label{eb2}
 H(G_{\epsilon}(t)) + \int_0^t R(G_{\e}(s)) \, \ds \leq H(G_{\epsilon}^{\text{in}}) \leq 
 C^{\text{in}} \epsilon^2 \,.
\end{equation}
From the definition of the relative entropy $H$ (see \eqref {enfunc} and \eqref{hg}), the terms that involve $H$ quantify how close $G_{\e}$ and $G_{\epsilon}^{\text{in}}$ are to the absolute equilibrium $G=1\,.$  Whereas, the terms involving the entropy production $R$ (dissipation rate specified in \eqref{endr}) can be regarded as formally measuring the proximity to equilibrium of the scaled collision integrand 
\begin{equation}\label{cint}
 q_{\epsilon} = \frac{1}{\e} (G'_{\epsilon} G'_{\epsilon*} -  G_{\epsilon} G_{\epsilon*})\,.
\end{equation}
As in \cite{BGL, Levermore}, the relative entropy and entropy production (dissipation rate) can be written as
\begin{align}\label{edrr}
 H(G_{\epsilon}) = \int_{\mathbb{T}^3} \langle h(\epsilon g_{\epsilon}) \rangle \, \dx \,, \qquad
 R(G_{\epsilon}) = \int_{\mathbb{T}^3} \frac{1}{4} \llangle[\bigg] 
 r \left(\frac{\epsilon q_{\epsilon}}{G_{\e} \,  G_{\e*}} \right)
 G_{\e} \,  G_{\e*} \rrangle[\bigg] \, \dx \,,
\end{align}
where $h$ and $r$ are the convex functions
\begin{align}\label{hr}
 h(z) = (1+z) \text{log} (1 + z) - z\,, \quad r(z) = z \text{log} (1+z)\,, \quad \tu{defined over } z > -1\,.
\end{align}
Since $h(z) = O(z^2)$ and $r(z) = O(z^2)$ as $z \to 0$, one clearly observes that $H(G_{\epsilon})$ 
and $R(G_{\epsilon})$ asymptotically act much like $L^2$ norms of $g_{\epsilon}$ and 
$q_{\epsilon}\,,$ respectively, as $\e \to 0\,.$  Using this notice, the following statement holds as a consequence only of the fact that the fluctuations satisfy the entropy bounds \eqref{eb2}.  It is a mere reformulation of Proposition 1 (\textit{the infinitesimal fluctuations lemma}) from \cite{Levermore} (derived from Proposition 3.1 and 3.4 of \cite{BGL}, which is followed by \cite[Lemma 5.1]{DL16} and \cite[Lemma 8.1]{4Diogo} also). 

\begin{proposition}\label{p1} 
Let $G_{\e} \equiv G_{\epsilon} (t,x,v) \geq 0$ be a sequence satisfying \eqref{gconditions} that 
\[G_{\epsilon}  \in C \left([0,\infty)\,;w\tu{-}L^1(M \dx \, \dv)\right) 
\cap L^{\infty}\left([0,\infty), \dt\,; L^1\left((1+|v|^2) M \dx \, \dv \right)\right)\,,\] 
where $G_{\e}$ fulfills the entropy inequality and bound (\ref{eb2}), with $G^{\textnormal{in}}_{\epsilon} 
= G_{\epsilon} (0)\,.$  Let $g_{\e}$ and $g^{\tu{in}}_{\e} = g_{\e} (0)$ be the corresponding 
fluctuations (\ref{ge}), and $q_{\epsilon}$ be the scaled collision integrands \eqref{cint} corresponding to $g_{\e}\,,$ with the chosen normalization $N_{\e}$ in \eqref{normalization}.  Then,
\begin{enumerate}[a)]
\item the sequence $g_{\e}(t)$ is relatively compact in $w\textnormal{-}L^1((1+|v|^2)M \dx \, \dv)$ for each $t \geq 0\,;$ 

\item if $g(t)$ is the $w\tu{-}L^1((1+|v|^2)M \dx \, \dv)$ limit of any converging subsequence of $g_{\e}(t)$ as $\epsilon \to 0\,,$ then $g(t) \in L^2(M \dx \, \dv)\,;$

\item the sequence $g_{\epsilon}$ is bounded in $L^{\infty}(\dt \,; L^1((1+|v|^2)M \dx \, \dv))$ and relatively compact in $w\tu{-}L^1_{\textnormal{loc}}(\dt \,; w\textnormal{-}L^1((1+|v|^2)M \dx \, \dv))\,;$ 

\item if $g$ is the $w\tu{-}L^1_{\textnormal{loc}}(\dt \,; w\textnormal{-}L^1((1+|v|^2)M \dx \, \dv))$ limit of any convergent subsequence of $g_{\e}\,,$ then $g \in L^{\infty}(\dt \,; L^2(M \dx \, \dv))\,;$

\item the sequence $q_{\epsilon}/N_{\epsilon}$ is relatively compact in $w\textnormal{-}L^1_{\textnormal{loc}}
 (\dt \,; w\tu{-}L^1((1+|v|^2)\dx \, \dmu))\,;$ 

\item if $q$ is the $w\textnormal{-}L^1_{\textnormal{loc}}
 (\dt \,; w\tu{-}L^1((1+|v|^2)\dx \, \dmu))$ limit of any converging subsequence of $q_{\epsilon}/N_{\epsilon}\,,$ then $q \in L^2(\dt \, \dx \, \dmu)\,;$
 
\item for almost every $t \geq 0$, any of these $g$ and $q$ satisfy the inequalities 
\begin{align}
& \int_{\mathbb{T}^3} \frac{1}{2} \langle g^2(t)\rangle \, \dx \leq \liminf\limits_{\epsilon \to 0} \int_{\mathbb{T}^3} \left\langle \frac{1}{\epsilon^2} h(\epsilon g_{\epsilon}(t)) \right\rangle \, \dx \leq  C^{\tu{in}} \,, \label{gie1} \\
& \int_0^t \int_{\mathbb{T}^3} \frac{1}{4} \llangle q^2 \rrangle \, \dx \, \ds \leq \liminf\limits_{\epsilon \to 0} \int_0^t \int_{\mathbb{T}^3} \frac{1}{4} \llangle[\bigg] \frac{1}{\epsilon^2} \, r\left(\frac{\epsilon q_{\epsilon}}{G_{\epsilon} G_{\epsilon*}} \right) G_{\epsilon} G_{\epsilon*} 
  \rrangle[\bigg] \, \dx \, \ds\,,  \label{qie1} \\
& \int_{\mathbb{T}^3} \frac{1}{2} \langle g^2(t) \rangle \, \dx + \int_0^t \int_{\mathbb{T}^3} 
\frac{1}{4} \llangle  q^2 \rrangle  \, \dx \, \ds \leq \liminf\limits_{\epsilon \to 0} \int_{\mathbb{T}^3} \left\langle \frac{1}{\epsilon^2} h(\epsilon g_{\epsilon}^{\textnormal{in}}) \right\rangle \, \dx \leq C^{\tu{in}} \,.  \label{gqie1} 
 \end{align}
\end{enumerate}
\end{proposition}

\begin{proof}
The proofs can be found in \cite[Propositions 3.1 and 3.4]{BGL}, especially, for \eqref{gie1} and \eqref{gqie1}.  The summaries of these proofs are presented in \cite[Lemmas 5.1 and 5.2]{DL16} and \cite[Lemmas 8.1 and 8.2]{4Diogo}, where the weak relative compactness of a sequence is the result of the boundedness, equi-integrability, and tightness in velocity of the sequence, in the manner of the Dunford-Pettis criterion \cite[p.~412]{Royden}.  Assertion \eqref{qie1} and its proof are summarized in \cite[Lemma~8.2]{4Diogo}, whose details can be found in \cite[Proposition 3.4 on pp.~702--703]{BGL}.   
\end{proof}

This Proposition \ref{p1} does not include the fact that the sequence $g_{\epsilon}$ will eventually describe 
fluctuations of the number density in the Boltzmann equation; and the proof is mainly based on 
the convexity of the integrands in the entropy inequality \eqref{gei2} 
as well as the entropy bound \eqref{eb2}, for the weak compactness.
\vspace{12pt}

We now revisit the notion of ``\textit{entropic convergence}'', which was first presented in \cite{BGL} and then in \cite{Levermore}.
\begin{definition} 
A family of fluctuations $g_{\epsilon}$ is said 
to \textbf{converge entropically} of order $\epsilon$ to $g \in L^2(M \dx \, \dv)$ (written as $ g_{\e}  \to g$ \textbf{entropically}) if and only if
\begin{align}
 \lim_{\e \to 0} g_{\epsilon} & = g \text{ in }  w\textnormal{-}L^1(M \dx \, \dv)\,, \label{ec1} \\
 \lim_{\e \to 0} \int_{\mathbb{T}^3} \left\langle \displaystyle \frac{1}{\epsilon^{2}} h(\epsilon g_{\epsilon})\right\rangle \, \dx 
& =
\int_{\mathbb{T}^3} \! \frac{1}{2}\langle g^2\rangle \, \dx\,. \label{ec2}
\end{align}
Using \eqref{edrr}, the condition (\ref{ec2}) is equivalent to
\begin{equation}\label{ec2b}
 \lim_{\e \to 0} \ \frac{1}{\epsilon^2} H(G_{\epsilon}) = \int_{\mathbb{T}^3} \! \frac{1}{2}\langle g^2\rangle \, \dx\,.
\end{equation}
\end{definition}

\begin{remark}\label{en2strong}
Throughout this paper, all entropic convergence will be assumed to be ``of order $\epsilon$''. Entropic convergence is a strong concept
because it implies that the sequence $g_{\e}$ converges to $g$ strongly in $L^1((1+|v|^2)M \dx \, \dv)\,,$ as demonstrated in Proposition 4.11 of \cite{BGL}.  That is, the notion of entropic convergence is stronger than that of the strong $L^1((1+|v|^2)M \dx \, \dv)$ topology. 
\end{remark}

Applying this concept of entropic convergence to Proposition \ref{p1}, the following 
sharpening of inequality (\ref{gqie1}) arrives, as a reformulation of Proposition 2 (\textit{the dissipation inequality corollary}) from \cite{Levermore} by a direct calculation \cite[p.~739]{BGL}.

\begin{proposition}\label{dic} 
Let $G_{\e} \equiv G_{\epsilon} (t,x,v) \geq 0$ be a sequence satisfying \eqref{gconditions} that
\[G_{\e} \in C \left([0,\infty)\,;w\tu{-}L^1(M \dx \, \dv)\right) 
\cap L^{\infty}\left([0,\infty), \dt\,; L^1\left((1+|v|^2) M \dx \, \dv \right)\right)\,.\] 
In particular, $G_{\epsilon}$ meets the entropy inequality and bound (\ref{eb2}).  We assume that $G^{\textnormal{in}}_{\epsilon} 
= G_{\epsilon} (0)$ has fluctuations $g_{\epsilon}^{\text{in}} = g_{\epsilon}(0)$ which converge entropically to some 
 $g^{\tu{in}} \in L^2(M \dx \, \dv)\,.$  Let $g_{\e}$ and $q_{\e}$ be the corresponding 
fluctuations (\ref{ge}) and scaled collision integrands \eqref{cint}; and let
 $g$ and $q$ be the corresponding weak limits.  Then,
\begin{equation}\label{edin}
\int_{\mathbb{T}^3} \frac{1}{2} \langle g^2(t) \rangle \, \dx + \int_0^t \int_{\mathbb{T}^3} 
  \frac{1}{4} \llangle  q^2 \rrangle \, \dx \, \ds \leq
  \int_{\mathbb{T}^3} \frac{1}{2} \langle (g^{\tu{in}})^2 \rangle \, \dx \,.   
 \end{equation}
\end{proposition}

\begin{remark}
Thanks to (\ref{ec2b}), the assumption that the initial fluctuations 
$g_{\epsilon}^{\textnormal{in}}$ converge entropically to some $g^{\textnormal{in}}$ in $L^2(M \dx \, \dv)$ indicates 
that these fluctuations meet the entropy bound (\ref{eb}).
\end{remark}

\bigskip
For removing the normalization as well as linearizing the collision integrand in the limit, the required technical 
estimations are obtained from the entropy inequality and bound \eqref{eb2}.  While the earlier propositions depend mostly on fundamental properties of convexity, the estimations 
in the following proposition rely mainly on the particular form of the relative entropy $H$ \eqref{enfunc}.  This proposition is contained within Corollary 3.2 and 
Proposition 3.3 in \cite{BGL} and is a reformulation of Proposition 3 (\textit{the bounds of nonlinear terms lemma}) from \cite{Levermore}. 

\begin{proposition}\label{bntl} 
Let $G_{\e} \equiv G_{\epsilon} (t,x,v) \geq 0$ be a sequence satisfying \eqref{gconditions} that 
\[G_{\epsilon}  \in C \left([0,\infty)\,;w\tu{-}L^1(M \dx \, \dv)\right) 
\cap L^{\infty}\left([0,\infty), \dt\,; L^1\left((1+|v|^2) M \dx \, \dv \right)\right)\,.\] 
More specifically, $G_{\e}$ fulfills the entropy inequality and bound (\ref{eb2}), with $G^{\textnormal{in}}_{\epsilon} 
= G_{\epsilon} (0)\,.$  Let $g_{\e}$ and $g^{\tu{in}}_{\e} = g_{\e} (0)$ be the corresponding 
fluctuations (\ref{ge}), where $g_{\e}$ has the normalized form $\gamma_{\e}$ defined in \eqref{gamma}, with the chosen normalization $N_{\e}$ in \eqref{normalization}.  Then,
\begin{enumerate}[a)]
\item if the sequence $g_{\epsilon}$ converges to $g$ in $w\textnormal{-}L^1_{\textnormal{loc}} (\dt \,; w\tu{-}L^1((1+|v|^2)M \dx \, \dv))\,,$ then the sequence $\gamma_{\epsilon}$ converges to $g$ in $w\textnormal{-}L^1_{\textnormal{loc}}(\dt \,; w\textnormal{-}L^1((1+|v|^2)M \dx \, \dv))\,;$

\item the sequence $g_{\epsilon}^2/N_{\epsilon}$ is bounded in $L^{\infty}(\dt\,; L^1(M \dx \, \dv))\,,$ and the bound is given by
\begin{equation}\label{bnt}
   \int_{\mathbb{T}^3} \left\langle \frac{g^2_{\epsilon}}{N_{\epsilon}} \right\rangle (t) \, \dx 
   \leq 2 C^{\textnormal{in}} \quad \forall t \geq 0\,;
\end{equation}

\item as $\epsilon \to 0\,,$
\begin{align}
|v|^2 \frac{g_{\epsilon}^2}{N_{\epsilon}} = O\left(\textnormal{log}\left(\frac{1}{\epsilon|\textnormal{log}(\epsilon)|}\right)\right) \quad \textnormal{in }
L^{\infty}\left(\dt\,;L^1(M \dx \, \dv)\right)\,.
\end{align}
\end{enumerate}
\end{proposition}

\bigskip
The statement $c$ above plays a crucial role in the following proposition, as a reformulation of Proposition 4 (\textit{the limiting collision integrand theorem}) from \cite{Levermore}, 
to control the quadratic terms 
of the scaled Boltzmann collision integrands \eqref{cint} in order to derive their linearized limit, as follows.  The proof can be found in \cite[pp.~1242--1243]{Levermore}.

\begin{proposition}\label{lcit}
Let $G_{\e} \equiv G_{\epsilon} (t,x,v) \geq 0$ be a sequence 
satisfying \eqref{gconditions} that 
\[G_{\epsilon}  \in C \left([0,\infty)\,;w\tu{-}L^1(M \dx \, \dv)\right) 
\cap L^{\infty}\left([0,\infty), \dt\,; L^1\left((1+|v|^2) M \dx \, \dv \right)\right)\,,\]
where $G_{\e}$ obeys the entropy inequality and bound (\ref{eb2}), with $G^{\textnormal{in}}_{\epsilon} 
= G_{\epsilon} (0)\,.$  Let $g_{\e}$ and $g^{\tu{in}}_{\e} = g_{\e} (0)$ be the corresponding 
fluctuations (\ref{ge}), and $q_{\epsilon}$ be the scaled collision integrands \eqref{cint} corresponding to $g_{\e}\,,$ with the chosen normalization $N_{\e}$ in \eqref{normalization}.  If the sequence $g_{\epsilon}$ converges to $g$ in $w\textnormal{-}L^1_{\textnormal{loc}}(\dt \,; w\textnormal{-}
 L^1((1+|v|^2)M \dx \, \dv))$ then
 \begin{align}\label{lcie}
  \frac{q_{\epsilon}}{N_{\epsilon}} \longrightarrow q = g' + g'_{*} - g - g_{*} \quad \textnormal{in }
  w\textnormal{-}L^1_{\textnormal{loc}}\left(\dt \,; w\textnormal{-}L^1(\dx \, \dmu)\right)\,.
 \end{align}
\end{proposition}

\bigskip
The preceding propositions apply only one fact about DiPerna--Lions solutions, which is 
the global entropy inequality (\ref{gei2}).  For a full proof, see \cite{Levermore}.  
The fact that they are also weak solutions of the renormalized 
Boltzmann equation (\ref{gammab}) leads to both the limiting dynamics, 
and subsequently a strong notion 
of convergence as in our following new proposition, which is a development of Proposition 5 (\textit{the limiting Boltzmann equation theorem}) from \cite{Levermore}, to include our considered external force term satisfying new condition $F \in L^1_{\tu{loc}}(\dt \,; L^2(M \dx\, \dv))\,.$

\begin{proposition}\label{lbet}
Let $F$ be a given force field  
\[F(t,x,v) \in L^1_{\textnormal{loc}}(\dt \, \dx \, \dv)\] 
such that 
\begin{equation}\label{Fconds1}
  \ddiv_v (F) = 0\,, \ F \cdot v = 0\,, \ \tu{and } 
  F \in L^1_{\tu{loc}}\left(\dt \,; W^{1,1}_{\tu{loc}}(\dx \, \dv)\right) \cap L^1_{\tu{loc}}(\dt \,; L^2(M \dx\, \dv)) \,.
 \end{equation}
Let $G_{\epsilon}$ be a family of the DiPerna--Lions renormalized solutions of the Boltzmann equation \eqref{vb2}--\eqref{vb20} (equivalently, weak solutions of the renormalized Boltzmann initial-value problem \eqref{rbe}--\eqref{rbe0}) satisfying the condition \eqref{gconditions} that
\[G_{\e} \in C \left([0,\infty)\,;w\tu{-}L^1(M \dx \, \dv)\right) 
\cap L^{\infty}\left([0,\infty), \dt\,; L^1\left((1+|v|^2) M \dx \, \dv \right)\right)\,.\]
Specifically, $G_{\e}$ fulfills the entropy inequality and bound (\ref{eb2}), having initial data $G_{\e}^{\tu{in}} = G_{\e} (0)$ that meets the entropy bound 
 (\ref{eb}).  Let $g_{\epsilon}$ and $g_{\e}^{\tu{in}} = g_{\e} (0)$ be the corresponding fluctuations (\ref{ge}).  Then,
 \begin{enumerate}[a)]
  \item the sequence $g_{\e}$ is relatively compact in 
  $C([0,\infty)\,; w\textnormal{-}L^1((1+|v|^2) M \dx \, \dv))\,;$

\item any convergent subsequence of $g_{\epsilon}$ has a limit $g \in C([0,\infty)\,; L^2(M \dx \, \dv))$ which is the unique solution of the initial-value problem
\begin{align}
\partial_t g + v \cdot \nabla_x g + F \cdot \nabla_v g + Lg = 0\,, \label{nlg} \\
g(0) = g^{\textnormal{in}} \equiv w\textnormal{-}L^1 \textnormal{-} \lim_{\epsilon \to 0} g_{\epsilon}(0)\,, \label{nlg0}
\end{align}
in which the linearized collision operator $L$ is defined in \eqref{lg}.
\end{enumerate}
\end{proposition}

\begin{proof}
The proof uses the strategy in \cite{Levermore, GL}.  First, we note that the global Maxwellian equilibrium $M \equiv M(v) >0$ defined in \eqref{m} is smooth and bounded in \eqref{Mbounds}. Thus, the condition $F \in L^1_{\tu{loc}}(\dt \,; W^{1,1}_{\tu{loc}}(\dx \, \dv))$ in \eqref{Fconds1} is equivalent to 
\begin{equation}\label{Fconds3}
F \in L^1_{\tu{loc}}\left(\dt \,; W^{1,1}_{\tu{loc}}(M \dx \, \dv)\right)\,.
\end{equation}

Recall the scaled collision integrand $q_{\e}$ defined in \eqref{cint} and the collision operator $Q(G,G)$ defined by \eqref{bco2}.  Consider the renormalized Boltzmann equation 
(\ref{gammab}), which is expressed as follows:
\begin{align}\label{ngamma}
 \partial_t \gamma_{\epsilon} + (v \cdot \nabla_x + F \cdot \nabla_v) \, \gamma_{\epsilon} 
 = \frac{1}{\epsilon} \frac{Q(G_{\epsilon}, G_{\epsilon})}{N_{\epsilon}} 
 =\int_{\mathbb{T}^3} \int_{\mathbb{S}^2} \frac{q_{\epsilon}}{N_{\epsilon}} \, b(v-v_{*}, \sigma) \, 
 \dsm \, M_{*} \, \dv_{*}\,.
\end{align}
Now, using \eqref{wsol} for $\gamma_{\e} = \dfrac{1}{\e} \beta(G_{\e})$ in \eqref{gamma}, the form \eqref{ngamma} means that with any $\varphi (x,v) \in W^{1,\infty}\left(M \dv\,;C^1(\mathbb{T}^3)\right)$ 
and every 
$ 0 \leq t_1 < t_2 < \infty\,,$ 
\begin{align}\label{nwle}
\begin{split}
\int_{\mathbb{T}^3} \langle (\gamma_{\epsilon}(t_2) - \gamma_{\epsilon}(t_1)) \, \varphi \rangle \, \dx 
- \int_{t_1}^{t_2} \int_{\mathbb{T}^3} \langle \gamma_{\epsilon}(v \cdot \nabla_x &+ F \cdot \nabla_v) \,  
\varphi \rangle \, \dx \, \dt \\
&=\int_{t_1}^{t_2} \int_{\mathbb{T}^3} \llangle[\bigg] \frac{q_{\epsilon}}{N_{\epsilon}} 
\varphi \rrangle[\bigg] \, \dx \, \dt \,.
\end{split}
\end{align}
Letting $\dd z= \frac{1}{3} \epsilon g_{\epsilon}$ in the elementary inequality 
\[\left(\text{log}(1+z)\right)^2 \leq \frac{z^2}{1+z} \quad \forall z > -1\,,\]
one obtains from the chosen normalization \eqref{normalization} that
\begin{equation}\label{gamineq}
\gamma_{\epsilon}^2 \leq \frac{g_{\epsilon}^2}{N_{\epsilon}}\,.
\end{equation}
The nonlinear bound (\ref{bnt}) in Proposition \ref{bntl} then demonstrates from \eqref{gamineq} that the sequence 
$\gamma_{\epsilon}$ is bounded in $L^{\infty}(\dt\,; L^2(M \dx \, \dv))\,,$ with
\begin{equation}\label{g2bound}
\int_{\mathbb{T}^3} \langle \gamma_{\epsilon}^2(t) \rangle \, \dx \leq 2C^{\text{in}} \quad \forall 
t \geq 0\,.
\end{equation}

Thanks to \cite[p.~365]{GL}, since the sequence $\gamma_{\epsilon}$ is bounded in $L^{\infty}(\dt\,; L^2(M \dx \, \dv))$ and $v 
\in L^2(M \dv)$ (indeed, $\int_{\mathbb{R}^3} |v|^2 M \dv < \infty$ by the form of $M$ defined in \eqref{m} and the normalization \eqref{norms}), it follows that the sequence $\langle v \gamma_{\epsilon} \rangle$ is relatively compact in 
$w\text{-}L^1_{\text{loc}}(\dt\,;w\text{-}L^2(\dx))$ ($\subset w\text{-}L^1_{\text{loc}}(\dt\,;w\text{-}L^1(\dx))$).  This can be proved by recalling the primary ideas from the proofs of \cite[Proposition 3.1]{BGL}, \cite[Lemma 5.1]{DL16}, and \cite[Lemma 8.1]{4Diogo}, or by using the steps below. 

With $F \in L^1_{\tu{loc}}(\dt \,; L^2(M \dx\, \dv))\,,$ we will show that the sequence 
$\langle F \gamma_{\e} \rangle$ is uniformly bounded and equi-integrable, so relatively compact in $w\text{-}L^1_{\text{loc}}(\dt\,;w\text{-}L^1(\dx))\,.$  Indeed, for any compact set $J = [t_1, t_2] \subset [0,\infty)\,,$ for almost everywhere $(t,x) \in J \times \mathbb{T}^3$ (with $| \mathbb{T}^3| =1$ from \eqref{norms}), by H{\"o}lder's inequality, we have
\[\int_{\mathbb{T}^3} |\langle F \gamma_{\e} \rangle(t,x)| \, \dx \leq \|F(t)\|_{L^2(M\dx \, \dv)} \ \|\gamma_{\e} (t)\|_{L^2(M\dx \, \dv)}\,.\]
Hence, utilizing the hypothesis $F \in L^1_{\tu{loc}}(\dt \,; L^2(M \dx\, \dv))\,,$ it holds that
\[\int_{J} \int_{\mathbb{T}^3} |\langle F \gamma_{\e} \rangle| \, \dx \, \dt \leq \sup_{\e > 0} \| \gamma_{\e} \|_{L^{\infty}(\dt\,; L^2(M\dx \, \dv))} \int_{J} \|F(t) \|_{L^2(M\dx \, \dv))} \dt < \infty\,.\]
Thus, the family $\{ \langle F \gamma_{\e} \rangle \}$ is uniformly bounded in $L^1(J \, ; L^1(\dx))\,.$ Next, we let 
\[C = \sup_{\e > 0} \| \gamma_{\e} \|_{L^{\infty}(\dt\,; L^2(M\dx \, \dv))} \,.\]
For any measurable set $A \subset J\,,$ the previous estimations yield
\begin{equation}\label{et}
\int_A \int_{\mathbb{T}^3} | \langle F \gamma_{\e} \rangle | \dx \, \dt \leq C \int_A \| F(t) \|_{L^2(M\dx \, \dv))} \, \dt \,.
\end{equation}
Since $F \in L^1_{\tu{loc}}(\dt \,; L^2(M \dx\, \dv))$ by the hypothesis \eqref{Fconds1}, it follows that the integral on the right-hand side of \eqref{et} has absolute continuity property \cite[Corollary 3.6 on p.~89]{folland}.  Hence, when $|A| \to 0\,,$ the right-hand side of \eqref{et} tends to 0, and thus the left-hand side of \eqref{et} tends to 0, uniformly in $\e\,.$  Therefore, the sequence $\{ \langle F \gamma_{\e} \rangle \}$ is uniformly integrable (also called equi-integrable) on $J \times \mathbb{T}^3$ (see \cite[Definition 5.2.2 on p.~104]{albiac_equiv}). Overall, according to the Dunford-Pettis criterion \cite[p.~412]{Royden}, we infer that the sequence $\langle F \gamma_{\epsilon} \rangle$ is weakly relatively compact in $L^1_{\text{loc}}(\dt\,; L^1(\dx))$ and hence relatively compact in $w\text{-}L^1_{\text{loc}}(\dt\,;w\text{-}L^1(\dx))\,.$ 

Therefore, in (\ref{nwle}), from left to right, the relative compactness of $\langle v \gamma_{\e} \rangle \,,$ $ \langle F \gamma_{\e} \rangle \,,$ 
and $q_{\epsilon} / N_{\epsilon}$ (see Proposition \ref{p1}) in $w\text{-}L^1$ (defined by \eqref{wf}) implies that each of the corresponding integrands in \eqref{nwle}, that is, $\langle \gamma_{\e} v \cdot \nabla_x \varphi \rangle\,,$ $\langle\gamma_{\e} F \cdot \nabla_v \varphi \rangle\,,$ and $\llangle (q_{\e}/N_{\e}) 
\varphi \rrangle$ is uniformly integrable on $[t_1,t_2] \times \mathbb{T}^3\,,$ which indicates that the map 
$t \mapsto \int_{\mathbb{T}^3} \langle \gamma_{\epsilon} (t) \varphi \rangle \, \dx$ is equicontinuous in $C([t_1,t_2])$ (defined in \eqref{wEt}) for each $\varphi$ (see also \cite[p.~1244]{Levermore}).
From \eqref{phisp}, note that $\varphi (x,v) \in W^{1,\infty}\left(M \dv\,;C^1(\mathbb{T}^3)\right) \subset L^{\infty}(M \dx \, \dv) \cong (L^1(M \dx \, \dv))^*$ and $\gamma_{\epsilon}(t) \in L^2(M \dx \, \dv) \subset L^1(M \dx \, \dv)$ (as $\dnu = M\dx \, \dv$ is a non-negative finite measure on $\mathbb{T}^3 \times \mathbb{R}^3$ of $(x,v)$ due to \eqref{norms}). It follows that the sequence $\gamma_{\epsilon}$, and thus $g_{\epsilon}$ is equicontinuous in 
$C([0,\infty)\,; w\text{-}L^1(M \dx \, \dv))$ (defined in \eqref{cYwE}).  The pointwise compactness of $g_{\epsilon}(t)$ is provided 
by Proposition \ref{p1}, hence the first claim holds by the Arzel\`{a}-Ascoli theorem \cite[p.~111]{brezis}.

\vspace{12pt}
For the second assertion, one passes to any convergent subsequence and take limits in (\ref{nwle}) as 
$\epsilon \to 0$ to reach 
\begin{align}\label{lwle}
 \begin{split}
\int_{\mathbb{T}^3} \langle g(t_2) \, \varphi \rangle \, \dx &- \int_{\mathbb{T}^3} \langle g(t_1) \, 
 \varphi \rangle dx - \int_{t_1}^{t_2} \int_{\mathbb{T}^3} \langle g(v \cdot \nabla_x + F \cdot 
 \nabla_v) \, \varphi \rangle \, \dx \, \dt \\
 &= \int_{t_1}^{t_2} \int_{\mathbb{T}^3} \llangle  q \, \varphi 
 \rrangle \, \dx \, \dt \,,
 \end{split}
\end{align}
where we have used Proposition \ref{bntl}, and $q = g' + g'_{*} - g - g_{*}$ by Proposition \ref{lcit}.  Substituting $q$ into the formula 
of $L$ in (\ref{lg}), one can demonstrate through \eqref{lwle} that $g$ is a weak solution of \eqref{nlg}--\eqref{nlg0}.  However then, $g$ must 
hence be the 
unique semigroup solution in $C([0,\infty)\,;L^2(M \dx \, \dv))\,.$
\end{proof}

\bigskip

Now, it comes that with a careful selection of initial data $G_{\e}^{\tu{in}}\,,$ any $L^2$ 
solution of the linearized Boltzmann initial-value problem \eqref{lbe}--\eqref{lbe0} 
can be obtained to uniquely describe 
the limiting fluctuation $g\,.$  In particular, notice that for every $g^{\text{in}} \in 
L^2(M \dx \, \dv)\,,$ it is always possible to build a non-negative sequence $G_{\epsilon}^{\text{in}} = 1 + 
\epsilon g_{\epsilon}^{\text{in}}$ such that its fluctuations $g_{\e}^{\tu{in}}$ 
converge entropically to $g^{\tu{in}}\,,$ for example by choosing $g_{\e}^{\tu{in}} 
= \tu{max}\{g^{\tu{in}}, -1/\epsilon\}$.  With this preparation of the initial data $G_{\epsilon}
^{\text{in}}\,,$ it holds in the following proposition that 
any approximating sequence of DiPerna--Lions fluctuations will even converge 
entropically for all positive times.

More specifically, from \cite{Levermore} by Levermore for the case without external force, 
we preserve the main result about the linearized limit, by developing Proposition 6 (\textit{the strong linearized limit theorem}) from \cite{Levermore}, to include our considered force term together with its new condition $F \in L^1_{\tu{loc}}(\dt \,; L^2(M \dx\, \dv))\,.$

\begin{proposition}[\textbf{Main Theorem}]\label{sllt}
Let $b(z, \sigma)$ be a locally integrable collision kernel satisfying the DiPerna--Lions assumption (\ref{dlassump}). Let $F$ be a given force field 
\[F(t,x,v) \in L^1_{\textnormal{loc}}(\dt \, \dx \, \dv)\] 
 such that 
\begin{equation}\label{Fconds2}
 \ddiv_v (F) = 0\,, \ F \cdot v = 0\,, \ \tu{and } 
  F \in L^1_{\tu{loc}}\left(\dt \,; W^{1,1}_{\tu{loc}}(\dx \, \dv)\right) \cap L^1_{\tu{loc}}(\dt \,; L^2(M \dx\, \dv)) \,.
 \end{equation}
 Given any initial data $g^{\textnormal{in}} \in L^2(M \dx \, \dv)\,,$ let 
 \[g \in C\left([0, \infty)\,; L^2(M \dx \, \dv)\right)\] 
 be the unique solution of the initial-value problem
\begin{align}\label{lbe2}
  \partial_t g + v \cdot \nabla_x g + F \cdot \nabla_v g + Lg = 0\,, 
  \qquad g(0) = g^{\textnormal{in}}\,,
 \end{align}
where the linearized collision operator $L$ is defined in \eqref{lg}. 
 For any sequence $G_{\epsilon}^{\textnormal{in}} = 1 + \epsilon g_{\epsilon}^{\textnormal{in}} \geq 0$ 
 such that 
$g_{\e}^{\tu{in}} \to g^{\tu{in}}$ entropically, let $G_{\epsilon}$ be any sequence of  renormalized solutions of the Boltzmann equation \eqref{vb2}--\eqref{vb20} (equivalently, weak solutions of the renormalized Boltzmann initial-value problem \eqref{rbe}--\eqref{rbe0}) satisfying \eqref{gconditions} that
\[G_{\e} \in C \left([0,\infty)\,;w\tu{-}L^1(M \dx \, \dv)\right) 
\cap L^{\infty}\left([0,\infty), \dt\,; L^1\left((1+|v|^2) M \dx \, \dv \right)\right)\,.\]
In particular, $G_{\e}$ fulfills the entropy inequality and bound (\ref{eb2}), with initial data $G_{\epsilon}^{\textnormal{in}} = G_{\e} (0)$ that meets the entropy bound (\ref{eb}). Let $g_{\epsilon}$ and $q_{\epsilon}$ be the corresponding fluctuations (\ref{ge}) and scaled collision integrands (\ref{cint}).  Then,
\begin{enumerate}[a)]
\item the sequence $g_{\epsilon}$ converges to $g$ in $C([0,\infty)\,;w\textnormal{-}L^1((1+|v|^2)M \dx \, \dv))$;

\item the sequence $g_{\epsilon}(t)$ converges entropically to $g(t)$ for each $t > 0$;
  
\item the sequence $q_{\epsilon}/N_{\epsilon}$ converges to $g' + g'_* - g - g_*$ in $L^1_{\textnormal{loc}}(\dt\,;L^1((1+|v|^2) \dx \, \dmu))$\,.
\end{enumerate}
\end{proposition}

\begin{proof}
The proof heavily relies on the new Proposition \ref{lbet}, which we rigorously proved above.  The rest of the proof follows \cite{Levermore} (without external force) and is presented here for completeness.  We first recall \eqref{Fconds3}.  

The proof for assertion 
$a$ is derived in the following manner.  Because $g_{\e}^{\tu{in}} \to g^{\tu{in}}$ entropically, by Remark \ref{en2strong}, $g_{\e}^{\tu{in}} \to g^{\tu{in}}$ in $w\tu{-}L^1((1+|v|^2)M \dx \, \dv)\,.$ Hence, by Proposition \ref{lbet}, every convergent subsequence of $g_{\e}$ must have the same limit, say $g\,,$ which is uniquely obtained from \eqref{lbe2}.  This concludes the first claim $a\,.$    

The proof for assertion 
$b$ follows from Proposition 
\ref{p1}, the notion of entropic convergence, the squeezing argument, 
the dissipation equality (\ref{disse}), and the entropy 
inequality part of (\ref{eb2}).  In particular, we begin with the following observations. Proposition \ref{p1} and its proof include that for any $t>0\,,$
\begin{align}\label{gie}
\begin{split}
& \int_{\mathbb{T}^3} \frac{1}{2} \langle g^2(t)\rangle \, \dx \leq \liminf\limits_{\epsilon \to 0} \int_{\mathbb{T}^3} \left\langle \frac{1}{\epsilon^2} h(\epsilon g_{\epsilon}(t)) \right\rangle \, \dx  \,,\\
& \int_0^t \int_{\mathbb{T}^3} \frac{1}{4} \llangle q^2 \rrangle \, \dx \, \ds \leq \liminf\limits_{\epsilon \to 0} \int_0^t \int_{\mathbb{T}^3} \frac{1}{4} \llangle[\bigg] \frac{1}{\epsilon^2} \, r\left(\frac{\epsilon q_{\epsilon}}{G_{\epsilon} G_{\epsilon*}} \right) G_{\epsilon} G_{\epsilon*} \rrangle[\bigg] \, \dx \, \ds\,.  
\end{split}
\end{align}
As $g_{\e}^{\tu{in}} \to g^{\tu{in}}$ entropically, \eqref{ec2} leads to $\displaystyle \lim_{\e \to 0} g_{\epsilon}  = g \text{ in }  w\textnormal{-}L^1(M \dx \, \dv)$ and
\begin{equation}\label{ec2in}
\lim_{\e \to 0} \int_{\mathbb{T}^3} \left\langle \displaystyle \frac{1}{\epsilon^{2}} h(\e g_{\e}^{\tu{in}})\right\rangle \, \dx 
=
\int_{\mathbb{T}^3} \! \frac{1}{2}\langle (g^{\tu{in}})^2\rangle \, \dx\,.
\end{equation}
Because $g$ is the solution of \eqref{lbe2}, it fulfills the dissipation equality \eqref{disse}, which makes \eqref{ec2in} become
\begin{equation}\label{ec2ds}
\lim_{\e \to 0} \int_{\mathbb{T}^3} \left\langle \displaystyle \frac{1}{\epsilon^{2}} h(\e g_{\e}^{\tu{in}})\right\rangle \, \dx 
=
\int_{\mathbb{T}^3} \frac{1}{2} \langle g^2(t) \rangle \, \dx +
\int_0^t \int_{\mathbb{T}^3} \frac{1}{4} \llangle  q^2 \rrangle \, \dx \, \ds \,.
\end{equation}
Writing the entropy inequality of \eqref{eb2} with the expressions \eqref{edrr}, we get
\begin{align}\label{eb2h}
\begin{split}
\int_{\mathbb{T}^3} \left\langle \frac{1}{\epsilon^2} h(\epsilon g_{\epsilon}(t)) \right\rangle \, \dx + \int_0^t \int_{\mathbb{T}^3} \frac{1}{4} \llangle[\bigg] \frac{1}{\epsilon^2} \, r\left(\frac{\epsilon q_{\epsilon}}{G_{\epsilon} G_{\epsilon*}} \right)  & G_{\epsilon} G_{\epsilon*} \rrangle[\bigg] \, \dx \, \ds \\
&\leq \int_{\mathbb{T}^3} \left\langle \displaystyle \frac{1}{\epsilon^{2}} h(\e g_{\e}^{\tu{in}})\right\rangle \, \dx\,.
\end{split}
\end{align}
Combining this inequality \eqref{eb2h} with \eqref{gie} and \eqref{ec2ds} indicates the existence of the following limits:
\begin{align}
\lim_{\epsilon \to 0} \int_{\mathbb{T}^3} \left\langle \frac{1}{\epsilon^2} h(\epsilon g_{\epsilon}(t)) \right\rangle \, \dx  
&= 
\int_{\mathbb{T}^3} \frac{1}{2} \langle g^2(t)\rangle \, \dx
\,, \label{last1}\\
\lim_{\epsilon \to 0} \int_0^t \int_{\mathbb{T}^3} \frac{1}{4} \llangle[\bigg] 
\frac{1}{\epsilon^2} \, r\left(\frac{\epsilon q_{\epsilon}}{G_{\epsilon} G_{\epsilon*}} \right) G_{\epsilon} G_{\epsilon*} \rrangle[\bigg] \, \dx \, \ds 
&=
\int_0^t \int_{\mathbb{T}^3} \frac{1}{4} \llangle q^2 \rrangle \, \dx \, \ds\,. \label{last2}
\end{align}
The expression \eqref{last1} means the entropic convergence claimed in $b\,.$  Note by \cite[Proposition 4.11]{BGL} that the entropic convergence in this assertion 
$b$ also implies the strong convergence of $g_{\e}(t)$ to $g(t)$ in $L^1((1+|v|^2)M \dx \, \dv)\,,$ as desired.

Last, the proof for assertion $c$ follows from \eqref{last2} and the final argument 
in the proof of Theorem 6.2 in \cite{BGL}.
\end{proof}

\begin{remark}
The proof of Proposition \ref{sllt} does not base on the velocity averaging theory, which was employed 
in achieving the basic global existence result of renormalized solutions \cite{LionsCauchy} and in examining the relative 
compactness  \cite{Lions1, Lions3} 
as well as in exploring fluid 
dynamical limits \cite{BGL, GL04, L, DL16}.  It is the regularity 
of the limiting linear dynamics as imposed in its dissipation \textit{equality} \eqref{disse} together 
with the idea of entropic convergence that derives the strong convergence above \cite{Levermore}.
\end{remark}

\section{Discussion}\label{discuss}
The result on strong linearized limit (Proposition \ref{sllt}) can hold with more general nonzero external forces $F$ as long as these forces satisfy at least both conditions $F \cdot v=0$ and $\div_v F =0$ and the regularity assumptions in \eqref{FcondsL2}. 

We have considered solely the Boltzmann equation and have not studied coupled systems, such as \cite{jabin_drag} and \cite{vasseur_drag} with external drag forces.  The notion of entropic convergence, however, can be applied to any kinetic equation possessing a convex entropy (thanks to Levermore), and with external forces as in this paper.  

Regarding open discussions, it is an expectation that giving the linearized Boltzmann equation with external force, its fluid dynamical limits \cite{Campini} can be established in the same spirit of our paper and \cite{Levermore}, deriving some kind of equations involving the external force.

This paper handles short-range interactions only.  It has not yet investigated, for kinetic equations with 
long-range interactions such as Coulombian 
interactions in plasmas, the linearization approximation, and the hydrodynamical 
transition toward some models of fluid mechanics (for instance, \cite{mai2024hd}).  On this trend, some questions arisen from 
the Vlasov--Maxwell--Boltzmann 
system have been nicely addressed by Ars\'{e}nio and Saint-Raymond 
\cite{DL16}, and several problems are still open.

\section{Conclusions}\label{conclude}
We have extended Levermore's results \cite{Levermore} regarding entropic convergence and the linearized limit for the Boltzmann equation (without external force) to the situations where external force is presented.  In particular, starting with an external force introduced in \cite[Section 4.1]{DL16} by Ars\'{e}nio and Saint-Raymond, we find new conditions on the external force, and carefully prove the preserving result in \cite{Levermore} with respect to the validity of the linearization approximation: provided any $L^2$ initial data for the linearized (about an equilibrium $M$) Boltzmann equation with external force, it can be demonstrated that any sequence of 
DiPerna--Lions 
renormalized solutions of the Boltzmann equation have fluctuations (about $M$) that converge entropically 
(and thus strongly in $L^1$) to the solution of the linearized Boltzmann equation 
for any positive time, given that its initial fluctuations about $M$ converge entropically to the $L^2$ initial data. 

In the future, even with different settings, it would be interesting to investigate an analog of entropic convergence notion in studying topics involving entropy, such as $H$-solutions \cite{Hsolution} or entropy decrease \cite{tadmor25entropydecrease}.  
  
\begin{appendix}
\section{The notion regarding spaces}\label{appendix}
This paper makes use of a number of fundamental topological linear spaces.  Some of our concepts regarding these spaces are conventional, whereas others are less common.  These spaces are explained below, along with our notation for them.  For a full description of them, we refer to \cite{BGL} and some standard 
references therein, for instance, see \cite{dunford_tvs}.  

Let $E$ be any normed linear space; $\|\cdot\|_E$ stands for its norm, and $E^{*}$ represents its dual space.  
Denote by $\langle \cdot \text{ }, \text{ } \cdot \rangle$ the natural bilinear form (or natural pairing, or canonical bilinear form) relating $E^{*}$ and $E\,.$  
The notion $w\text{-} E$ will be used to indicate the space $E$ equipped with its weak topology, that is 
the coarsest topology on $E$ for that any 
of the linear forms
\begin{equation}\label{wEt}
u \mapsto \langle v , u \rangle_{E^{*},E} \qquad (\text{for } v \in E^{*})
\end{equation}
is continuous. 

Let $X$ be a locally compact topological space and $E$ a normed linear space.  The standard notion $C(X\,;w\text{-}E)$ will be used to indicate the space of continuous functions from $X$ to $w\text{-}E$, 
that is the set of function $u$ for which each of the linear forms
\begin{equation}\label{cYwE}
x \mapsto \langle v , u(x) \rangle_{E^{*},E} \qquad \text{is in } C(X) \qquad \text{for each } 
v \in E^{*}\,. 
\end{equation}
On such spaces $C(X\,;w\text{-}E)\,,$ we note that the Arzel\`{a}-Ascoli theorem is valid.

Let $(Y, \Sigma, \nu)$ be a measure space and $E$ a Banach space (complete normed linear space or complete normed vector space).  For every $1 \leq p \leq 
\infty$, we will use the abbreviated notion $L^p(\dnu \,;E)$ for the Bochner space 
$L^p((Y, \Sigma,\dnu)\,;E)$, containing all functions $u: Y \to E$ such that the corresponding norm 
is finite:
\[\|u\|_{L^p(\dnu \,;E)} := \left(\int_Y \|u(y) \|_E^p \, \dnu(y) \right)^{1/p} < + \infty 
\quad \forall \quad 1\leq p < \infty\,,\]
\[\|u\|_{L^{\infty}(\dnu\,; E)}: = \text{ess sup}_{y \in Y} \|u(y)\|_E = 
\text{min}\{\alpha: \nu(\{y: \| u(y)\|_E > \alpha\}) =0\} < + \infty\,.\]
Note that for $1 \leq p \leq \infty\,,$ the Bochner spaces $L^p(\dnu \,;E)$ are Banach spaces.

Also, we have
\[C([0,\infty);E)= \{ u \ | \ u: [0,\infty) \to  E \tu{ continuous} \}\,,\]
which is also a Banach space with the norm $\| u \|:= \displaystyle \max_{t \in [0,\infty)} \| u(t) \|_E\,.$  It can be considered as a closed subspace of the Banach space
\begin{align*}
L^{\infty}([0,\infty);E) = \{ [u] \ | \ u: [0,\infty) \to E &\tu{ is Bochner measurable and }\\ 
& \tu{ess sup}_{t \in [0, \infty)} \| u(t) \|_E < \infty \tu{ holds}\}\,.
\end{align*}
Furthermore,
\[C\left([0,\infty);L^2(Y, \Sigma,\dnu)\right) =\{ u \ | \ u:[0,\infty) \to L^2(Y, \Sigma,\dnu) \tu{ continuous} \}\,.\]

We will use $L^p(\dnu)$ to denote the same space $L^p(\dnu \,;E)$ 
whenever $E$ is a power of $\mathbb{R}$.  

For $1 \leq p < \infty$, the dual space of $L^p(\dnu\,;E)$ is $L^{p^{*}}(\dnu \,;E^{*})$, where 
$p^{*} = p/(p-1)$.  In this paper, only $p=1,2, \infty$ arise.

When $Y$ is locally compact, $\dnu$ is a Borel measure, and $E$ is any Banach space, we will denote by $L^p_{\text{loc}} 
(\dnu \,; E)$ (or $L^p_{\text{loc}}(\dnu)$), the space specified by the family of seminorms
\[u \mapsto \left( \int_K \|u(y) \|^p_E \, \dnu(y)\right)^{1/p} \qquad \text{for compact } K \subset Y\,, \quad 1\leq p < \infty\,.\]

For every $1 \leq p < \infty\,,$ we will use the notation $w\text{-}L^p(\dnu \,;w\text{-}E)$ 
(or $w\text{-}L^p(\dnu)$) to denote the space $L^p(\dnu \,; E)$ equipped with its weak topology, which is 
the coarsest topology on $L^p((Y, \Sigma, \dnu)\,;E)$ for that any of the linear forms
\begin{align}\label{wf}
u \mapsto \int_Y \varphi(y) \,  \langle v , u(y) \rangle_{E^{*},E} \, \dnu(y) \quad (\text{for } 
v \in E^{*} \quad \text{and } \varphi \in L^{p^{*}}(\dnu))
\end{align}
is continuous.  Hereafter, we will utilize the prefix $w$- to express that a given space is equipped
with its weak topology. 

When $Y$ is locally compact, for any $1 \leq p < \infty\,,$ we denote by $w\text{-}L^p_{\text{loc}}(\dnu \,; w\text{-}E)$ 
(or $w\text{-}L^p_{\text{loc}}(\dnu)$) the space $L^p_{\text{loc}}(\dnu \,; E)$ endowed with its weak 
topology; that is, the coarsest for which all the linear form (\ref{wf}) are continuous, where 
the functions $u$ are restricted to possess compact support in $Y\,.$
\end{appendix}


\section*{Acknowledgments}
The author would like to thank Professors Claude Bardos, Russel Caflisch, Fran\c{c}ois Golse, David Levermore, and Laure 
Saint-Raymond very much for really helpful conversations which have also motivated this work (since 2016, UT Austin).  She especially thanks Prof.~Bardos and Prof.~Levermore for encouraging comments, Prof.~Saint-Raymond for her kindheartedness and invaluable advices on the 
first draft of the paper, and Prof.~Caflisch for inspiration and suggestions about the external forces (at IPAM, UCLA 2016, the QKP2023, IMS, NUS, Singapore, and Frontiers in Computational Mathematics Conference 2025, UT Austin, US).  She is indebted to the sponsoring by Prof.~Yalchin Efendiev and the Institute for Scientific Computation, at Texas
A\&M University, College Station, US (TAMU).  Her respectful gratitude to Prof. Edriss Titi for motivating discussions and important encouragement (TAMU 2025).  She genuinely thanks Prof.~Minh-Binh Tran for friendly meetings and open-minded conversations.  Also, members of TAMU Math Department and the Department Head Prof.~Peter Howard 
offer a warm hospitality for the author during completing this paper (2023, 2025).  Her appreciations to the colleagues she has had the chances to meet at conferences, to people, institutions and places hosted her visits 
during 2015--2016, having the University of Minnesota, and UT Austin, where the work was originated from. Her many thanks to the Ki-Net Conferences, US, 2015--2016, and the IMS program Multiscale Analysis and Methods for Quantum and Kinetic Problems, NUS, Singapore (QKP2023) with generous support and good colleagues.  The work was done partially while the author was participating in the program QKP2023 of the Institute for Mathematical Sciences, National University of Singapore, in 2023.

\bibliographystyle{plain} 
\bibliography{ecvb}
\end{document}